\documentclass[10pt]{amsart}
\usepackage{amsmath,amssymb}
\usepackage[cp1251]{inputenc}
\usepackage{graphicx}

\newtheorem{thm}{Theorem}[section]
\newtheorem{prop}[thm]{Proposition}

\newtheorem{defi}[thm]{Definition}
\newtheorem{lem}[thm]{Lemma}

\theoremstyle{definition}

\newcommand{\x}{\mathbb{X}}
\newcommand{\T}{\mathbb{T}}                     %
\newcommand{\V}{{\mathcal V}}                    %
\newcommand{\op}[1]{\operatorname{#1}}             % operatorname
\newcommand{\mps}[3]{#1:#2\rightarrow #3}           % mapping
\newcommand{\matop}[3]{\mathop{#1}\limits_{#2}^{#3}} % mat oper.
\newcommand{\pam}[4]{#1\hookleftarrow #2\mathop{\rightarrow}
\limits^{#3}#4}                                 % partial mapping
\newcommand{\ii}{^{-1}}                          % ^{-1}
\newcommand{\mean}{\rightleftharpoons}          % mean
\newcommand{\ANE}{\operatorname{ANE}}          % ANE
\newcommand{\diam}{\operatorname{diam}}         % diameter
\newcommand{\cov}{\operatorname{cov}}           % open covers
\newcommand{\Id}{\operatorname{Id}}             %
\begin{document}
\title[On extending actions of groups]
{On extending actions of groups}
\author[Sergei M. Ageev]{Sergei M. Ageev}
\address[Sergei M. Ageev]{Faculty of Mathematics and Mechanics,
 Belarus State University, Independence Avenue 4, Minsk, Belarus 220050}
\email{ageev\_serhei@yahoo.com}
\author[Du\v{s}an Repov\v{s}]{Du\v{s}an Repov\v{s}}
 \address[Du\v{s}an Repov\v{s}]{Institute of Mathematics, Physics and Mechanics, and Faculty of Education, University of Ljubljana,  P. O. B. 2964 , Ljubljana, Slovenia 1000}
\email{dusan.repovs@guest.arnes.si}
\begin{abstract} We study a
problem proposed by E.V. Shchepin, concerning extensions of actions of compact transformation groups, under various
assumptions. We present several applications of our methods which we develop in this paper.
\end{abstract}
\thanks{The authors were supported in part by Slovenian Research Agency grants P1-0292-0101-04 and J1-9643-0101.}
\maketitle

%%%%%%%%%%%%%%%%%%%%%%%%%%%%%%%%%%%%%%%%%%%%%%%%%%%%%%%
%%%%%%%%%%%%%%%%%%%%%%%%%%%%%%%%%%%%%%%%%%%\input{Intr}

\section{Introduction}
 The diagram ${\mathcal D}\mean\{\Bbb
X\mathop{\rightarrow}\limits^{p}X\mathop{\hookrightarrow}\limits^{i}
Y\}$ is called  {\it admissible } if  $p\colon\Bbb X\rightarrow X$
is  an orbit projection and $i$ is a topological embedding of the
orbit space  $X$ into a hereditarily paracompact space $Y$. We say
that the {\it problem of extending the action is solvable for the
admissible diagram $\mathcal D$} if  there exists an equivariant
embedding $j:\Bbb X\hookrightarrow\Bbb Y$ into $G$-space $\Bbb Y$
(which is called a {\it solution} of the problem of extending the
action for a given diagram) covering $i$, that is,  the embedding
$\tilde{j}:X\hookrightarrow p(\Bbb Y)$ of orbit spaces induced by
$j$ coincides with $i$ (in particular, $p(\Bbb Y)=Y$).

We say that the {\it problem of extending the action} (denoted
briefly by {\it PEA }) is {\it solvable} if  there exists a solution
of the PEA  for each admissible diagram $\mathcal D$, i.e.
the diagram $\mathcal D$ yields
a commutative 
diagram:
$$
\begin{array}{ccc}
\mathbb{X} & \stackrel{j}{\hookrightarrow} & \mathbb{Y}\\
p\downarrow && \downarrow p\\
X & \stackrel{i}{\hookrightarrow} & Y.
\end{array}
$$

The problem of extending  group actions  naturally splits
into
{\it the closed and dense parts} -- depending on
the type of the
embedding $i$. On the other hand, it is clear that the
simultaneous solvability of the closed and dense PEA implies the
solvability of the PEA  in general.

 In what follows we shall
 consider only compact groups  $G$. 
 Motivated by his own results, 
 E.V. Shchepin was
the first to state the closed problem of extending  group  actions.
  
  \begin{prop}\label{one+2} If  the closed problem of
extending the group action   is solvable  for all  {\it metric
admissible diagrams}, i.e. the diagrams in which
$\x$ and $Y$ are  metric (it is clear  that such a diagram is
admissible). Then the following holds:
 \begin{enumerate}
 \item[$(1)$]
 The orbit  space of any
$G$-$\op {A[N]E} (\mathcal M)$-space $\Bbb E$ with metrizable
orbits is an  $\op {A[N]E}({\mathcal M})$-space (here  the class
of all metric $G$-spaces is denoted by $\mathcal M$).
\end{enumerate}
 \end{prop}
 
 {\it Proof.} Let  $\pam{Z}{A}{\varphi}{E}$
be a partial map with $Z\in{\mathcal M}$. We denote 
the fiberwise product $A_\varphi\times_{p}\Bbb E$
by $\Bbb A$. Since  all
orbits of $\Bbb A$ are metrizable and $\Bbb A/G=A\in\mathcal M$,
$\Bbb A\in\mathcal M$. In view of $Z\in{\mathcal M}$, the
embedding $A\hookrightarrow Z $ is covered by a  closed
$G$-embedding $\Bbb A\hookrightarrow\Bbb Z\in{\mathcal M}$. 
Consider the partial $G$-map $\pam{\Bbb Z}{\Bbb A}{\varphi'}{\Bbb
E}$, where $\varphi'$ is parallel to $\varphi$. Since $\Bbb E\in
G$-$\op {A[N]E}({\mathcal M})$, the $G$-map $\varphi'$ can be
$G$-extended into $\Bbb Z$ [into a $G$-neighborhood]. Passing to
the orbit spaces we get the desired extension $\hat\varphi$ of
$\varphi$. $\ \ \square$

 \medskip
 It turns out that the converse fact is also true  \cite{Ag11,Ag1}.
 
 \begin{prop}\label{one+12}
 The validity of the condition (1) implies the solvability of the  closed PEA  for all
metric admissible diagrams.
 \end{prop}
 
 {\it Proof.} Let $\T_G$ be the countable product of the metric cone $\op{Con}
\T$ over a discrete union of all transitive spaces $G/H\in
G$-$\ANE$. By Theorem \ref{v.03.04+}
there exists an isovariant
map $\mps{f}{\Bbb X}{\T_G}$. It is clear that $\T_G\in
G$-$\op{AE}({\mathcal M})$. Hence it follows by (1) that
$T_G=\T_G/G\in \op{AE}({\mathcal M})$ and the partial map
$\pam{Y}{X}{\tilde{f}}{T_G}$ induced by $f$ has an extension
$\mps{\varphi}{Y}{T_G}$. It is easy to verify that the fiberwise
product $Y_{\varphi}\times_p\T_G$ is a desired $G$-space $\Bbb
Y\in{\mathcal M}$. $\ \ \square$

\bigskip The dense PEA arises
naturally in theory of equivariant compactifications and it was
first stated by L. Zambahidze and  Yu. M. Smirnov. By now it 
has been settled
for the class of metric spaces with an action of an
arbitrary
zero-dimensional compact group \cite{Ag2}. We remark without proof
that the dense PEA is intimately connected with the existence of several
invariant metrics on a given $G$-space.

 \begin{prop}\label{one+2+1}
The dense problem of extending the group action is solvable for all
{\it separable metric diagrams}, i.e. the diagrams in
which  $\x$ and $Y$ are separable metric spaces, if and only if the following holds:
 \begin{enumerate}
 \item[$(2)$]
For each separable metric $G$-space $\Bbb X$ and each compatible
metric $d$ on $X$ there exists a compatible  invariant metric
$\rho$  on $\Bbb X$ such that ${\mathcal F}_d={\mathcal
F}_{\tilde{\rho}}$ (here ${\mathcal F}_d$ is the set  of all
sequences fundamental with respect to $d$, and $\tilde{\rho}$ is
an induced metric on $X$).
\end{enumerate}
 \end{prop}
 
\medskip 

The main purpose of the present paper is to present a
direct
proof of the following theorem which gives the positive answer  to
the closed PEA under very general assumptions.

 \begin{thm}\label{one}
The closed problem of extending the group action is solvable for
all  admissible diagrams.
  \end{thm}
  
We remark that the closed PEA fails to be solvable in certain simple
situations. Let  $\Bbb X=(\op{Con}\Bbb Z_2)^{\omega_1}$ be the
semi-free $\Bbb Z_2$-space, and $i:X\hookrightarrow I^{\omega_1}$
an arbitrary embedding into Tihonoff cube of uncountable weight
$\omega_1$ (note that $I^{\omega_1}$ is not hereditarily
paracompact). It can be checked that the PEA  is unsolvable for
the diagram $\Bbb
X\mathop{\rightarrow}\limits^{p}X\mathop{\hookrightarrow}\limits^{i}
I^{\omega_1}$.

\medskip 

If  the acting group $G$ is non-metrizable, then
some orbits of the solution of the closed PEA  for metric
admissible diagram $\mathcal D$ given by Theorem \ref{one} can
be nonmetrizable. The  following result eliminates this defect.

  \begin{thm}\label{one*}
  If the $G$-space $\Bbb X$ in the closed
admissible diagram  $\Bbb
X\mathop{\rightarrow}\limits^{p}X\mathop{\hookrightarrow}\limits^{i}
Y$ is stratifiable, then there exists  a solution $s\colon\Bbb
X\hookrightarrow\Bbb Y$ of the PEA  for this  diagram such that
$\Bbb Y\subset Y\times \Bbb Z$ for some metric $G$-space $\Bbb Z$.
  \end{thm}
  
In particular,   
the direct proof of solvability
of the  closed PEA  in the class of metrizable  $G$-spaces, as
well as in the class of stratifiable $G$-spaces, follows.

\medskip 

Theorem \ref{one} and the method developed for its proof
admit various  modifications. We endow the set $\mathcal S$ of all
solutions of the closed PEA for an
admissible diagram ${\mathcal
D}=\{\Bbb
X\mathop{\rightarrow}\limits^{p}X\mathop{\hookrightarrow}\limits^{i}
Y\}$ by
the  following partial order:  $s_2\ge s_1$ where
$\{s_i\colon\Bbb X\hookrightarrow\Bbb Y_i\}\subset\mathcal S$ if
and only if there exists a $G$-map $\mps{h}{\Bbb Y_1}{\Bbb Y_2}$
such that $h\circ s_1=s_2$ and $\tilde h=\op{Id}_{Y}$.

 \begin{thm}\label{one**+}{\rm (On majority)}
Each  finite set  of solutions $\{s_i\}_{i=1}^{n}\subset\mathcal
S$ has a majority, i.e. there exists  $s\in\mathcal
S$ such that $s\ge s_i$ for each  $i\le n$.
  \end{thm}
  
 \begin{thm}\label{one**}{\rm (On special solutions)}
For each  closed admissible diagram $\mathcal D$ there exists a
solution  $s:\Bbb X\hookrightarrow\Bbb Y$ of the PEA for $\mathcal
D$ such that for each $y\in\Bbb Y$, $\{ x\in\Bbb X\mid
(G_y)\ge(G_x)\}\not=\emptyset$.
  \end{thm}
  
 We shall omit the proofs of Theorem  \ref{one**+}
and Theorem   \ref{one**} because of their similarity with
the proof of Theorem \ref{one}. In their turn, Theorems \ref{one*},\
\ref{one**+} and \ref{one**} imply several results on the
equivariant extensor theory (which are, as can be shown,
equivalent to the corresponding facts  on extending of action).
The reader should reconstruct the proofs on his own.

 \begin{thm}\label{one12}
 If the  stratifiable $G$-space $\Bbb X$ is an equivariant absolute
extensor  for the class of stratifiable spaces, then its orbit
space $X$ is an absolute extensor  for the class of stratifiable
spaces.
  \end{thm}
  
 In \cite{AgRep} the proof of this
theorem was reduced to the Borsuk-Whitehead-Hanner theorem for
stratifiable spaces. Unfortunately, the status of the latter
theorem is still open in view of the gap  in \cite{Cau}. In
connection with establishing Theorem \ref{one12}, it would be
interesting to return to \cite{Cau} and to settle the
Borsuk-Whitehead-Hanner theorem for stratifiable spaces overcoming
the gap.

 \begin{thm}\label{one13}
If  the metric $G$-space $\Bbb X$ is an equivariant absolute
extensor, then for all
finitely many $G$-extensions
$\mps{h_i}{\Bbb Z}{\Bbb X}$ of the partial  $G$-map
$\pam{{\mathcal M}\ni\Bbb Z}{\Bbb A}{f}{\Bbb X}$ there exists  a
$G$-extension $\mps{h}{\Bbb Z}{\Bbb X}$ of $f$ such that
$(G_{h(z)})\ge (G_{h_i(z)})$ for any $z\in \Bbb Z$ and $i$.
  \end{thm}
  
   \begin{thm}\label{one14}
 Let  $G$ be a compact group,
$\x\in G$-$\ANE(\mathcal M)$ and ${\mathcal C}\subset\op{Orb}G$.
Then $\x^{\mathcal C}=\{x\mid (G_x)\ge (H)\ \text{for some}\
(H)\in {\mathcal C}\}\subset\x$ is a $G$-$\ANE(\mathcal M)$.
  \end{thm}
  
Earlier this theorem was proved by M. Murayama \cite{Mu}
for all
compact Lie group and one-element
collection ${\mathcal C}$.

 \medskip 
 
 In conclusion we present
 a result revealing the role of the set
of extensor points in the equivariant extensor theory. 
We say that
the $G$-embedding $\Bbb Y\hookrightarrow\Bbb X$ is {\it
equivariantly homotopically dense}, if  there exists  a $G$-homotopy
$\mps{F}{\Bbb X\times[0,1]}{\Bbb X}$ such that $F_0=\op{Id}$ and
\begin{enumerate}
 \item[$(\alpha)$] $F_t(\Bbb X)\subset\Bbb Y$ for each  $t>0$.
 \end{enumerate}
 
   \begin{thm} \label{v.03.04+3}{\rm (On equivariant homotopy density)}
 Let  $\Bbb X$ be a metric $G$-$\ANE$-space. Then
the subspace $\Bbb X_{\mathcal E}$ of all  extensor points is
equivariantly homotopically dense in $\Bbb X$.
  \end{thm}
  
   As a immediate consequence we get the  following result:
 \begin{thm} \label{one**++}
 Let  $\Bbb X$ be a metric $G$-$\ANE$-space. Then
each  $G$-subspace  $\Bbb Y,\Bbb X_{\mathcal E}\subset\Bbb
Y\subset\Bbb X$, is a $G$-$\ANE$.
  \end{thm}
  
%%%%%%%%%%%%%%%%%%%%%%%%%%%%%%%%%%%%%%%%%%%%%%%%%%%%%%%%%%%%%%%%%%%
%%%%%%%%%%%%%%%%%%%%%%%%%%%%%%%%%%%%%%%%%%%%%% \input Veden

\section{Preliminaries}

Let $G$ be a compact group. An {\it action of $G$ on a space $X$}
is a homomorphism $T:G\to\op{Aut} X$ of $G$ into the group
$\op{Aut} X$ of all autohomeomorphisms of $X$ such that the map
$G\times X\to X$ given by $(g,x)\mapsto T(g)(x)=g\cdot x$ is
continuous. A space $X$ with a fixed action of $G$ is called a
{\it $G$-space}.

For any point $x\in X$ the {\it isotropy subgroup} of $x$, or the
{\it stabilizer} of $x$, is defined as $ G_x=\{g\in G\mid g\cdot
x=x\}$ and the orbit of $x$ as $ G(x)=\{g\cdot x\mid g\in G\}$.
The space of all orbits is denoted by $X/G$ and the natural map
$p=p_{\Bbb X}:X\to X/G$, given by $p(x)=G(x)$, is called the {\it orbit
projection}. The orbit space $X/G$ is equipped with the quotient
topology induced by $p$. In what follows we shall
denote $G$-spaces
and their orbit spaces variously: $\Bbb X,\Bbb Y,\Bbb Z,\dots $
for $G$-spaces, and $X,Y,Z,\dots$ for their orbit spaces.

The map $\mps f{\Bbb X}{\Bbb Y}$ of $G$-spaces is called {\it
equivariant} or {\it $G$-map}, if $f(g\cdot x)=g\cdot f(x)$ for
all  $g\in G$ and $x\in\x$. Each $G$-map $\mps f{\Bbb X}{\Bbb Y}$
induces a map $\mps {\tilde f}{X}{Y}$ of orbit spaces by formula
$\tilde f(G(x))=G(f(x))$. Equivariant map $\mps f{\x}{\Bbb Y}$ is
said to be {\it isovariant}, if $G_x=G_{f(x)}$ for all $x\in\Bbb
X$.

  \begin{thm}{\rm (Equimorphism criterion  \cite{Br})}
\label{CriEquim} An isovariant map $\mps{f}{\x}{\Bbb Y}$ inducing
the homeomorphism $\mps{\tilde f}{X}{Y}$ of orbit spaces is an
equimorphism.
\end{thm}

 Observe
 that all $G$-spaces and $G$-maps generate a category
denoted by  $G$-$\op{TOP}$ or $\op{EQUIV}$, provided that no
confusion occurs. If  "$***$" is a well-known notion from
nonequivariant topology, then "$G$-$***$" or "$\op{Equiv}$-$***$"
means the corresponding equivariant analogue. See \cite{Br} for
more details on compact  transformation groups.

 The subset  $A\subset\Bbb X$ is called {\it invariant} or {\it
$G$-subset}, if   $G\cdot \Bbb A=\Bbb A$. For each  closed
subgroup $H$ of $G$ (briefly $H<G$) we define the  following sets:
$\x^H = \{x\in\x \mid H\cdot x = x\}$ (which is called {\it
$H$-fixed points set}); $\x^{(H)} = \bigcup\{\x_K \mid K<G \
\text{and}\ H'<K\ \text{for some subgroup}\ H'\ \text{ conjugated
with} \ H\} = G\cdot \x^H$; $\x_{(H)} = \bigcup\{\x_K \mid K<G \
\text{and}\ K \hbox{ is }\text{conjugated  with } H\} = G\cdot
\x_H$.

We endow the set $\op{Orb}(G)$ of  all  conjugate classes of
closed subgroups of $G$ with the  following partial order:
$(K)\preccurlyeq (H)\Leftarrow\Rightarrow K\subset g\ii \cdot
H\cdot g\ \text{for some }\ g\in G$. It is clear that
$\op{type}\Bbb X\mean\{(G_x)\mid x\in\Bbb X \}$ is a subset of
$\op{Orb}(G)$ (hereafter the sign $\mean$ is used for
the introduction of the new objects placed to the left of it). If
${\mathcal C}\subset\op{Orb}(G)$, then $\x^{\mathcal
C}\mean\{x\mid (G_x)\ge (H)\ \text{for some}\ (H)\in {\mathcal
C}\}\subset\x$.

 We now introduce several concepts related to extensions of
 $G$-maps partially defined in $G$-spaces from the class
$\mathcal K$ -- one of the increasing chain: the class of all
metric $G$-spaces $\subset$ the class of all stratifiable
$G$-spaces {the space $X$ is  {\it stratifiable} if there
exists a family $\{\mps{f_U}{X}{[0,1]}\mid U\subset X\ \text{is an
open subset }\}$ of continuous functions such that $f_U\ii(0,1]=U$
and $f_U\le f_V$ if and only if $U\subset V$) $\subset$ the class
$\mathcal P$ of all hereditarily paracompact $G$-spaces
(the space is called {\it hereditarily paracompact}  if each
its subspace is paracompact, i.e. is equivalent to the
paracompactness of each its open subspace). 

A space X is called
an {\it absolute neighbourhood extensor for class $\mathcal K$},
$\x\in\op G$-$\op{ANE}({\mathcal K})$, if each $G$-map
$\mps{\varphi}{\Bbb A}{\x}$ defined on a closed $G$-subset $\Bbb
A\subset \Bbb Z$ of $G$-space  $\Bbb Z\in \mathcal K$ and called
the {\it partial $G$-map} can be $G$-extended in a
$G$-neighborhood $\Bbb U\subset \Bbb Z$ of $\Bbb A$,
$\mps{\hat\varphi}{\Bbb U}{\x}, \hat\varphi\restriction_{\Bbb
A}=\varphi$ (we use the notation $f\restriction_A$ for
the restriction of the map $f$ to $A\subset X$ or we simply write
$f\restriction$ if the set $A$ in question is clear). 

If it is
possible to $G$-extend  $\varphi$ in $\Bbb U=\Bbb Z$, then $\x$ is
called an {\it equivariant absolute extensor for class $\mathcal
K$}, $\x\in\op G$-$\op{AE}({\mathcal K})$. If the acting group $G$
is trivial, then these notions are transformed into the notions of
absolute [neighborhood] extensors for class $\mathcal K$ --
$\op{A[N]E}({\mathcal K})$. Since we are mainly interested in
equivariant absolute [neighborhood] extensors for the class of
metrizable $G$-spaces, we will briefly denote them as
$G$-$\op{A[N]E}$. The following results are well-known: each
Banach $G$-space \cite{Matu}, each compactly convex $G$-subset
of locally convex complete vector $G$-space \cite{Ab},p.155, each
linear normed $G$-spaces (for acting compact Lie group $G$)
\cite{Mu},p.488, are $G$-$\op{AE}$-spaces.

 We will  heavily depend
 on the
the slice theorem \cite{Br} which is equivalent to the following
assertion: "If $G$ is a compact Lie group, then the transitive
space $G/H\in G$-$\op{ANE}$ for the class of regular $G$-spaces".

We recall the construction of the equivariant absolute extensor
for arbitrary compact group  $G$ go back to \cite{Ma}. Recall that
$G$ acts on the space $\x=\op{C}(G,Y)$ of all continuous maps
endowed with compact-open topology by formula $(g\cdot
f)(h)=f(g^{-1}\cdot h)$ where $f\in\op C(G,Y)$ and $g,h\in G$. If
$Y$ is metric, then $\x$ is also metric.

 \begin{thm}\label{v.1.4} If  metric
space $Y$ is  $\op{AE}$-space for class $\mathcal P$ of
paracompact spaces, then  $\op{C}(G,Y)$ is a $G$-$\op{AE}$-space
for $\mathcal P$.
\end{thm}

By \cite{Dk} each  Banach space $B$ is a $G$-$\op{AE}$ for class
of paracompact spaces. Hence it follows by Theorem \ref{v.1.4}
that $\op C(G,B)$ is $G$-$\op{AE}$ for $\mathcal P$.

 \begin{prop} \label{t5.9++} {\rm (Palais Metatheorem \cite{Pal})} 
 Let ${\mathcal P}(H)$
be a property  which depends on compact Lie group  $H$. Suppose
that ${\mathcal P}(H)$ is true, provided ${\mathcal P}(K)$ is true
for each compact Lie group $K$ isomorphic to a proper subgroup of
$H$. If ${\mathcal P}(H)$ is true for trivial group $H=\{e\}$,
then ${\mathcal P}(H)$ is true for all compact Lie groups $H$.
 \end{prop}
 
 By a {\it fiberwise product} of the spaces $C$ and $B$ with respect to
maps $g$ and $f$ we call the  subset   $\{(c,b)\mid
g(c)=f(b)\}\subset C\times B$ which is denoted by $C_g\times_fB$.
The projections $D=C_g\times_fB$  onto the factors $C$ and $B$ are
denoted by $\mps{\check{f}}{D}{C}$ and $\mps{\check{g}}{D}{A}$.
These maps $\check{f}$ and $\check{g}$ are called the {\it maps
parallel} to $f$ and $g$ respectively, and we write for brevity
$\check{f}\| f$ and $\check{g}\| g$. It should remarked that the map
$\mps{f\circ \check{g}= g\circ \check{f}}{D}{A}$ is the product of
$g$ and $f$ in the category {\bf TOP}$_A$ of all spaces over $A$.
The most important example of the fiberwise product in the theory
of compact transformation groups is supplied by isovariant maps.

 \begin{prop} \label{pull 3.1} Let
$\mps{h}{\Bbb Y}{\Bbb X}$ be an isovariant map,
$\mps{\tilde{h}}{Y}{X}$ the map of orbit spaces induced by $h$.
Then $\Bbb Y$ is the fiberwise product
$Y_{\tilde{h}}\times_{p_\Bbb X}\Bbb X$, moreover  $h\| \tilde{h}$
and $p_{\Bbb Y}\| p_{\Bbb X}$ where $p_{\Bbb Y},p_{\Bbb X}$ are
the orbit projections.
 \end{prop}

\bigskip

We denote by ${\mathcal D}$ the
commutative square diagram in the category $\op{EQUIV}$ with
closed $\Bbb A\subset\Bbb Z$:
 $$
\begin{array}{ccc}
\mathbb{A} & \matop{\longrightarrow}{}{\varphi} & \mathbb{X}\\
\cap && \downarrow f\\
\mathbb{Z}& \matop{\longrightarrow}{}{\psi} &\mathbb{Y}
\end{array}
$$
 We will say that the {\it $G$-map $\varphi$ is the partial
lifting of the $G$-map $\psi$ with respect to } $f$. The {\it
problem of extension of partial lifting for ${\mathcal D}$ is
solvable globally  $[$locally $]$} if there exists a $G$-map
$\mps{\hat\varphi}{\Bbb Z}{\Bbb X}$ $[\mps{\hat\varphi}{\Bbb
U}{\Bbb X}$ defined on a neighborhood $\Bbb U\subset \Bbb Z$ of
$\Bbb A]$ such that $\hat\varphi\restriction_{\Bbb A}=\varphi$ and
$f\circ\hat\varphi=\psi$
$[f\circ\hat\varphi=\psi\restriction_{\Bbb U}]$. We will say also
that $\mps{\hat\varphi}{\Bbb Z}{\Bbb X}$ $[\mps{\hat\varphi}{\Bbb
U}{\Bbb X}]$ is a global $[$local$]$ lifting of $\psi$.

  \begin{defi} \label{pull IsEx-10-1}
A morphism  $\mps{f}{\Bbb X}{\Bbb Y}$ in the category $\op{EQUIV}$
is called  {\it equivariantly soft $[$locally  equivariantly soft$]$} if
for each  commutative square diagram  ${\mathcal D}$ in the
category $\op{EQUIV}$ the problem of extending of partial lifting
is solvable  globally $[$locally $]$.
\end{defi}

 \subsection{$P$-orbit projection.}
Let  $P$ be a normal closed subgroups of $G$ (briefly
$P_\alpha\triangleleft G$) and $\mps{\pi}{G}{H\mean
G/P},\pi(g)=g\cdot P$, be a natural epimorphism. It is clear that
$x\sim x' \Leftrightarrow x'\in P\cdot x$ is a equivalent relation
on the $G$-space $\x$. Then the quotient space $\x/P$ defined by
this relation coincides with $\{P\cdot x\mid x\in\x\}$. It is
clear that $\Bbb Y=\x/P$ is an $H$-space: $(g\cdot P)\cdot(P\cdot
x)= P\cdot(g\cdot x)$. If $y=P\cdot x$, then the stabilizer $H_y$
coincides with $G_x\cdot P$.

The quotient map $\mps{f}{\x}{\x/P}$ is called  the {\it $P$-orbit
projection}. If   $P=G$, then  $f$ coincides with the orbit
projection $\mps{p}{\x}{\x/G}$. Since the composition of the
$P$-orbit projection $f$ and the $H$-orbit projection $\Bbb Y$ is
perfect, $f$ is  a perfect surjection and satisfies the following
properties:
\begin{enumerate}\setcounter{enumi}{0}
  \item $f(gx)=\pi(g)\cdot f(x)$ for all  $x\in\x$ and
 $g\in G$;
 \item $\pi(G_x)=H_{f(x)}$ for all  $x\in\x$; and
 \item If $f(x)=f(x')$, then $x$ and $x'$ belong to the same orbit.
 \end{enumerate}

\noindent The  following fact shows that these  properties
characterize $P$-orbit projections completely.

   \begin{prop}\label{Proj1}  Let
$\mps{\pi}{G}{H}$ be a epimorphism of compact groups,
$P=\op{Ker}\pi$. The perfect surjection $\mps{f}{\x}{\Bbb Y}$ from
the $G$-space $\x$ onto the $H$-space $\Bbb Y$ is the $P$-orbit
projection if and only if the properties $(1)-(3)$ for $f$ hold.
 \end{prop}
 
 {\it Proof.} We consider the $P$-orbit projection
$\mps{\varphi}{\x}{\x/P=\Bbb Z}$ and define the map
$\mps{\theta}{\Bbb Z}{\Bbb Y}$ by the formula  $\theta(P\cdot
x)=f(x)$. It follows by $(1)$ that $\theta$ is correctly defined
and it is equivariant. It follows by perfectness of $f$ that
$\theta$ is a perfect surjection, and therefore the induced map
$\tilde{\theta}$ of orbit space is also perfect and surjective.

 Let  $z=P\cdot x$ and $z'=P\cdot x'$. If  $\theta(z)=\theta(z')$ then  $f(x)=f(x')$. In view of
 $(3)$,  $x$ and $x'$ lie on the same $G$-orbit. Therefore  $z$ and $z'$ lie on the same
$G/P$-orbit. Hence $\tilde{\theta}$ is a homeomorphism.

 The map $\theta$ preserves the orbit type of  points as on the one hand,
$H_{\theta(z)}=H_{f(x)}$, and on the other hand, $H_{z}=G_x\cdot
P=\pi(G_x)$. By $(2)$ we have $H_{\theta(z)}=H_{z}$. Therefore
$\theta$ is an isovariant map inducing on the orbit spaces a
homeomorphism. By Theorem \ref{CriEquim} the map $\theta$ is an
equimorphism. $\ \square$

\bigskip For compact group $G$ we consider the {\it Lie series $\{
P_\alpha\triangleleft G\}$} of normal closed subgroups indexed by
ordinals $\alpha<\omega$ \cite{Ptr}. This means that
\begin{enumerate}\setcounter{enumi}{4}
 \item $P_1=G$; $P_\beta<P_{\alpha}$ for all
$\alpha<\beta$; $P_\alpha/P_{\alpha+1}$ is a compact Lie group for
all $\alpha<\omega$, $P_\alpha=\cap \{P_\alpha'\mid
 \alpha'<\alpha\}$ for each  limit ordinal $\alpha$, and also $\cap \{P_\alpha\mid
 \alpha<\omega\}=\{e\}$. \end{enumerate}

In this  case $G$ is the limit $\matop{\lim}{\leftarrow}{}\{
G/P_\alpha, \varphi^{\beta}_{\alpha}\} $ of the inverse system of
quotient groups $\{G/P_\alpha\}$ and  natural epimorphisms
$\mps{\varphi^{\beta}_{\alpha}}{G/P_{\beta}}{G/P_\alpha},\alpha<\beta$.
A more general fact holds.

   \begin{lem}\label{End1}
Let  $\mps{f^{\beta}_{\alpha}}{\Bbb X_{\beta}}{\Bbb X_\alpha}$ be
natural projections from $\Bbb X_{\beta}=\x/P_{\beta}$ into
$\break\Bbb X_{\alpha}=\x/P_{\alpha}$. Then  $f^{\beta}_{\alpha}$
is a $P_{\alpha}/P_{\beta}$-orbit projection and the map
$\mps{f}{\x}{\matop{\lim}{\leftarrow}{} \{ \x_\alpha,
f^{\beta}_{\alpha}\}}$, given by the formula $f(x)=\{P_\alpha\cdot
x\}$, is an equimorphism.
 \end{lem}
 
 The {\it proof} of Lemma \ref{End1} consists of a
 straightforward application
of the equimorphism criterium \ref{CriEquim}. The converse  to
Lemma \ref{End1} is also valid.

  \begin{lem}\label{End1+}
Let  $\{ P_\alpha\triangleleft G\}$ be the Lie series, $\
P_{\alpha}^{\beta},\alpha<\beta$, a kernel of the homomorphism
$\varphi^{\beta}_{\alpha}$, and $\mps{g^{\beta}_{\alpha}}{\Bbb
Z_{\beta}}{\Bbb Z_\alpha}$ a $P^{\beta}_{\alpha}$-orbit projection
with $g^{\beta}_{\alpha}\circ g^{\gamma}_{\beta}=
g^{\gamma}_{\alpha}$ for all $\alpha<\beta<\gamma$. Then  $\Bbb Z
=\matop{\lim}{\leftarrow}{} \{ \Bbb Z_\alpha,
g^{\beta}_{\alpha}\}$ is  a $G$-space and $\Bbb Z_\alpha=\Bbb
Z/P_\alpha$.
 \end{lem}
 
%%%%%%%%%%%%%%%%%%%%%%%%%%%%%%%%%%%%%%%%%%%%%%%%%%%%%%%%%%%%%%%%%%%
%%%%%%%%%%%%%%%%%%%%%%%%%%%%%%%%%%%%%%%%%%%%%%%%\input{Extgroup}

\section{Extensor subgroups}

  \begin{defi} \label{aslice2+7} A closed subgroup $H<G$ of
  a compact
group $G$ is called a ${\mathcal P}$-subgroup if  the transitive
space $G/H$ is finite-dimensional and locally  connected.
 \end{defi}
 
 L.S. Pontryagin \cite{Ptr} proved that  $H<G$ is a ${\mathcal P}$-subgroup
 if and only if one of the  following properties holds:
 
  \begin{enumerate}
   \item[$(1)$] There exists a normal subgroup $P\triangleleft G$
such that $P<H$ and $G/P$ is a compact Lie group;
 \item[$(2)$] $G/H$ is a topological manifold.
\end{enumerate}

 It is known \cite{Ptr} that each  compact group contains
 arbitrarily small normal ${\mathcal P}$-subgroups. Hence
the  following fact is valid:

 \begin{prop}\label{slice2-7} Let  $G$ be a compact group, $\Bbb X$ a
$G$-space. Then
  \begin{enumerate}\setcounter{enumi}{2}
   \item for each  neighborhood
${\mathcal O}(H)\subset G$ of the subgroup $H<G$ there exists a
${\mathcal P}$-subgroup $H'\triangleleft G$ such that $H\subset
H'\subset{\mathcal O}(H)$;
 \item for each  neighborhood ${\mathcal O}(x)$ of the point $x\in\Bbb X$ there exists
a normal $\mathcal P$-subgroup $P\triangleleft G$ such that the
$P$-orbit projection $\mps{p}{\Bbb X}{\Bbb X/P\mean \Bbb Y}$ has a
small inverse image of $y\mean P\cdot x$, $p^{-1}(y)\subset
{\mathcal O}(x)$. Moreover,  there exists a neighborhood
${\mathcal W}\subset \x/P$ of $y$ such that $p^{-1}({\mathcal W}
)\subset {\mathcal O}(x)$.
\end{enumerate}
 \end{prop}
 
 \noindent The property (3) easily implies that
  \begin{enumerate}\setcounter{enumi}{4}
 \item $G/H$ is metric if and only if  $H<G$ can be presented as intersection of
 countably many ${\mathcal P}$-subgrous.
\end{enumerate}

If  $H<G$ is a ${\mathcal P}$-subgroup then  it follows by $(1)$
and the slice theorem \cite{Br} that $G/H$ is $G$-$\ANE$. The
proof of the converse fact is based on the existence of a regular
$G$-space $\Bbb Z$ such that the stabilizers of all its points are
${\mathcal P}$-subgroups except a nowhere dense orbit $G(z)\cong_G
G/H$. Hence it follows that:
 \begin{enumerate}
 \item[$(6)$]  $H<G$ is a ${\mathcal P}$-subgroup if and only if   $G/H$.
\end{enumerate}
The equivalence $(6)$ expresses the main property of ${\mathcal
P}$-subgroups and it
simultaneously  justifies the alternative term
 -- {\it extensor subgroups}. We also draw  reader's attention to
 the  following result \cite{Sz}: "If the natural action of the compact group $G$
on $G/H\in\ANE$ is effective, then $G$ is a
Lie group", from
which it follows that
\begin{enumerate}
 \item[$(7)$]  $H<G$ is a ${\mathcal P}$-subgroup if and only if   $G/H$ is
an $\ANE$.
\end{enumerate}
   \begin{defi} \label{aslice2-5}
We call the conjugate class $(H)$  {\it extensor}, provided $H<G$
is an extensor subgroup. By $\Bbb X_{\mathcal E}$ we denote the
collection of all {\it extensor points of $\Bbb X$}, that is, all
points $x\in\Bbb X$ for which $G_x<G$ is an extensor subgroup.
 \end{defi}

 We say that the $G$-subspace  $\Bbb Y\subset\Bbb X$ is   {\it $G$-dense}, if
$\Bbb Y^H\subset\Bbb X^H$ is dense for each   subgroup $H<G$. In
\cite{AgU} it was shown that
\begin{enumerate}\setcounter{enumi}{7}
 \item  $\Bbb X_{\mathcal E}\subset\Bbb X$ is   $G$-dense if and
only if $\Bbb X$ is a $G$-$\ANE$ for the class of all metric
$G$-spaces with zero-dimensional orbit space;
 \item (equivariant Dugunji's
theorem) a linear normed $G$-space $\Bbb L$ is a $G$-$\op{AE}$ if
and only if $\Bbb L_{\mathcal E}\subset\Bbb L$ is $G$-dense.
\end{enumerate}
 There exists an example of a linear normed $G$-space $\Bbb L\not\in G$-$\op{AE}$
for which $\Bbb L_{\mathcal E}\subset\Bbb L$ is  dense (but not
$G$-dense).

\smallskip We list the basic properties of extensor subgroups. It can be
shown that
\begin{enumerate}\setcounter{enumi}{9}
 \item  the subgroup $H<G$ is  extensor if and only if $G$ admits an
orthogonal action on $\Bbb R^n$ such that $G/H$ is equimorphic to
the orbit of a point.
\end{enumerate}

It is clear  that each subgroup $H<G$ in compact Lie group $G$,
and also  each  clopen subgroup $H<G$ are extensor subgroups.

It easily follows by $(1)$ that the property to be an extensor
subgroup is inherited by any passage to the larger subgroup. It is
well-known that
a
compact group is a Lie group if and only if it
contains no small subgroups  \cite{Ptr}. Hence it follows that the
quotient group $G/(P_1\cap P_2),P_i\triangleleft G$, are Lie group
if and only if each $G/P_i$ is a Lie group.
  \begin{prop}\label{aslice2-6}
The intersection of finitely many extensor subgroups is an
extensor subgroup.
 \end{prop}
 This proposition cannot be improved: if the extensor subgroup $H<G$ is
an intersection of a family $\{H_\alpha<G\}$ of extensor
subgroups, then $H$ is an intersection of finitely many subgroups
$H_{\alpha_i}$.

Since each  closed subgroup of compact Lie group is again  compact
Lie group \cite{Ptr}, it easily follows that if  $H<G$ is an
extensor subgroup, then
\begin{enumerate}\setcounter{enumi}{9}
 \item $H\cap K<K$ is an extensor subgroup for each $K<G$; and
 \item  $H<K$ and $K<G$ are extensor subgroups for each $K<G,H<K<G$.
 \end{enumerate}
 The proof of the next
 fact follows  from the definition of extensor subgroup and
Hurewicz  theorem  on dimension \cite{AP}.
  \begin{prop} \label{aslice2-7}
Let  $K<L<G$ and $K<L$ be an extensor subgroup. Then $K<G$ is an
extensor subgroup if and only if  $L<G$ is an extensor subgroup.
 \end{prop}

 \medskip If   $G$ is a compact non-Lie group, then
by  $(7)$ $G\not\in\op{ANE}$. It is known that there exists a free
$G$-space $\Bbb X\in\op{ANE}$ \cite{AgIzv}. Hence each its
invariant open subset  is not homeomorphic to a product $G\times
U$ and therefore in this  case the slice theorem fails \cite{Vil}.
But if we weaken the requirement about the slice, then the  following is
valid. We say that a $G$-map $\mps{\alpha}{\Bbb X}{G/H},H<G$, is
{\it sliced} if $G/H\in G$-$\op{ANE}$ ($\equiv$ $\ H<G$ is an
extensor subgroup).

    \begin{thm} \label{slice2-6} {\rm (On approximate slice of
$G$-space  \cite{Ab},p.151,\cite{AgU,Ag4}))} Let a compact
group $G$ act on a $G$-space $\Bbb X$. Then for each  neighborhood
${\mathcal O}(x)$ of  $x\in\Bbb X$ there exists a  neighborhood
$\V=\V(e)$ of the unit $e\in G$, an extensor subgroup $K<G,G_x<K$,
and a slice map $\mps{\alpha}{\Bbb U}{G/K}$ where $\Bbb U$ is an
invariant neighborhood of $x$ such that $x\in\alpha\ii({\mathcal
V}\cdot [K])\subset{\mathcal O}(x)$.
 \end{thm}
 
 {\it Proof  of Theorem \ref{slice2-6}.}
 Let $P\triangleleft G$ be a normal $\mathcal P$-subgroup,
${\mathcal W}$ a neighborhood of $y=P\cdot x\in \Bbb Y\mean\Bbb
X/P$ taken from Proposition \ref{slice2-7}$(4)$. Since  $\Bbb Y$
is naturally endowed with action of compact Lie group $G'=G/P$,
and $G'_y=G_x\cdot [P]<G'$, there exists  by slice theorem a
neighborhood $\V(e)$ and a slice map $\mps{\alpha}{\Bbb
V}{G'/G'_y\cong G/(G_x\cdot P)}$ defined on a $G'$-neighborhood
$\Bbb V$ of the orbit $G'(y)$ such that
$\alpha\ii(\V(e)\cdot[G'_y])\subset {\mathcal W}$. We easily check
that $K\mean G_x\cdot P<G$ is an extensor subgroup and the
composition $\mps{\alpha\circ p}{p\ii(\Bbb V)}{G/K\in G$-$\ANE}$
is the desired slice map. $\ \ \square$

\medskip The  following result is known for compact Lie groups $G$
\cite{Dick},7.6.4.
 \begin{thm} \label{slMap2-1} Let   $H$ and $K$ be subgroups of a
compact group $G$ such that $H<K$ is an extensor subgroup. Then
the natural projection $\mps{p}{G/K}{G/H}$ is equivariantly
locally soft.
 \end{thm}

%%%%%%%%%%%%%%%%%%%%%%%%%%%%%%%%%%%%%%%%%%%%%%%%%%%%%%%%%%%%
%%%%%%%%%%%%%%%%%%%%%%%%%%%%%%%%%%%%%%%%\input{Prel}

 \section{Tube structure of orbit projections}
 
 We consider the epimorphism   $\mps{\pi}{G}{H}$ of compact groups with
the kernel $P$ being  Lie group. Let  $K<G$ be an extensor
subgroup  and $\pi(K)\mean L<H$. Since  $K<\pi\ii(L)<G$, it
follows by \ref{aslice2-7}$(10)$ that
\begin{enumerate}
 \item $K<\pi\ii L$  and $\pi\ii(L)<G$
are  extensor subgroups,
\end{enumerate}
 and hence
 \begin{enumerate}\setcounter{enumi}{1}
 \item $G/K\in G$-$\ANE$ and $G/\pi\ii(L)\cong H/L\in H$-$\ANE$.
\end{enumerate}
Let  $\mps{\kappa}{G/K}{H/L}$ be a composition of the natural
epimorphism $\mps{\alpha}{G/K}{G/\pi\ii L}$ and the isomorphism
$\mps{\beta}{G/\pi\ii L}{H/L}$. By  $(1)$ and Theorem
\ref{slMap2-1} we have the following
 \begin{enumerate}\setcounter{enumi}{2}
 \item The maps  $\mps{\alpha}{G/K}{G/\pi\ii
L}$ and $\mps{\kappa}{G/K}{H/L}$ are  equivariantly locally soft.
\end{enumerate}
Since  $\kappa(g\cdot [K])=\pi(g)\cdot [L]$ and   $\pi(K)=L$, it
follows that
 \begin{enumerate}\setcounter{enumi}{3}
 \item $\kappa(g\cdot [K])=[L]$ if and only if  $g\in P\cdot K$.
\end{enumerate}

 We say that the $P$-orbit projection $\mps{f}{\Bbb
X}{\Bbb Y}$ have a {\it  $\kappa$-tube structure generated by
slice maps} $\mps{\varphi}{\Bbb X}{G/K \in G$-$\ANE}$ and
$\mps{\psi}{\Bbb Y}{H/L\in H$-$\ANE}$ if  they close the following
{\it diagram $\mathcal A$} up to commutative one:
$\kappa\circ\varphi=\psi\circ f$.
  $$
\begin{array}{ccc}
\Bbb X & \stackrel{\varphi}{\rightarrow} &G/K\\
f\downarrow&&\downarrow \kappa\\ \Bbb
Y&\stackrel{\psi}{\rightarrow} &H/L
\end{array}
$$
 The $\kappa$-tube structure on $f$ is said to be {\it
nontrivial} if  $\kappa$ is not  bijection. It is equivalent to
$K\subsetneqq \pi\ii L$ or $P\setminus K\not=\emptyset$.

\begin{prop}\label{Proj2-+} Let
$\mps{f}{\Bbb X}{\Bbb Y}$ be a  $P$-orbit projection. If
$x\not\in\Bbb X^P$, then the restriction of $f$ on the orbit
$G(x)$ has a nontrivial tube structure.
\end{prop}
 {\it Proof of Proposition \ref{Proj2-+}.}
Since   $P\setminus G_x\not=\emptyset$ then  by  Proposition
\ref{slice2-7}$(3)$ there exists  an extensor  subgroup $K<G$ such
that $P\setminus K\not=\emptyset$ and   $G_x<K$. It is clear that
$G_{f(x)}$ coincides with $G_x\cdot P$. The  following subgroups
$G_x<K<K\cdot P$ and $G_x<G_{f(x)}<K\cdot P$ naturally generate
$G$-maps $\mps{\varphi}{G(x)}{G/K}$,
$\mps{\psi}{G(f(x))}{G/(K\cdot P)}$ and
$\mps{\kappa}{G/K}{G/(K\cdot P)}$, which in its turn generate the
nontrivial tube structure on the orbit $G(x)$. $\ \ \square$

\medskip We consider the diagram $\mathcal A$ and  denote by
$\Bbb S$ and $\Bbb T$ the $K$-space $\varphi\ii([K])$ and the
$L$-space $\psi\ii([L])$, correspondingly. It is clear that $\x$ is
$G$-homeomorphic to the twisted product $G\times_K \Bbb S$ and
$\Bbb Y\cong_H H\times_L\Bbb T$. The following facts show that if
$f$ has a nontrivial tube structure, then $f$ is generated by a
$Q$-orbit projection where $Q$ is a proper  subgroup of $G$.

  \begin{lem}\label{Proj2}
Let $f$ has a tube structure given by epimorphism $\kappa$. Then
$f(\Bbb S)=\Bbb T$ and  the map $\mps{f\restriction}{\Bbb S}{\Bbb
T}$ is a $Q$-orbit projection where  $Q$ is the kernel  of the
epimorphism $\mps{\pi'=\pi\restriction}{K}{L}$. If this tube
structure  is nontrivial, then  $Q$ is a proper subgroup of $P$.
\end{lem}
 {\it Proof.} Let  $y\in\Bbb T$, $x\in f\ii(y)$ and $\varphi(x)=g\cdot K$. Then
$[L]=\psi(y)=\psi(f(x))=\kappa(\varphi(x))=\kappa(g\cdot K)$. By
$(4)$, $g\in P\cdot K$. Now it is easily to check that $x'\mean
g\ii \cdot x\in\Bbb S$ and $f(x')=y$. Hence  $f(\Bbb S)=\Bbb T$.
It is clear  that the map  $f\restriction$ is perfect.

 Let us  verify the properties  $(1)-(3)$ of Proposition
\ref{Proj1} for restriction  $f\restriction$: the property $(1)$
and the rest of the properties  $(3)$  is performed obviously. The
property $(2)$ holds, as $K_s=G_s$ and $L_{f(s)}=H_{f(s)}$ for
$s\in\Bbb S$.

 If $\kappa$ is not a bijection, then $K\subsetneqq \pi\ii L$ and hence  $P\setminus
K\not=\emptyset$. But $Q=\op{Ker}\pi'$ coincides with $K\cap P$,
and  hence it is proper subgroup of $P$.  $\ \square$

 \medskip We now consider the converse situation:  there exist
an epimorphism  $\mps{\pi'}{K}{L}$ of compact groups and a
$Q$-orbit  projection $\mps{f'}{\Bbb S}{\Bbb T}$ where $\Bbb S$ is
a $K$-space, $\Bbb T$ is an $L$-space and $Q=\op{Ker}\pi'$. Let
$K<G$ and $L<H$ be extensor subgroups and $\mps{\pi}{G}{H}$ an
epimorphism  extending  $\pi'$. By $\mps{\kappa}{G/K}{H/L}$ we
denote the composition of natural epimorphisms
$\mps{\alpha}{G/K}{G/\pi\ii L}$ and $\mps{\beta}{G/\pi\ii
L}{H/L}$. Then the formula $f([g,s]_K)=[\pi(g),f'(s)]_L$ correctly
defines the map $\mps{f=\pi\times f'}{G\times_K \Bbb
S}{H\times_L\Bbb T}$. It is straightforwardly  checked  that the
perfectness of $f'$ implies the same property for $f$.
  \begin{lem}\label{Proj3} The map
$\mps{f}{G\times_K \Bbb S}{H\times_L\Bbb T}$ is a $P$-orbit
projection  for $P=\op{Ker}\pi$. In so doing the natural slice
maps $\mps{\varphi}{G\times_K \Bbb
S}{G/K},\mps{\psi}{H\times_L\Bbb T}{H/L}$ and also  epimorphism
$\mps{\kappa}{G/K}{H/L}$ set a tube  structure  on the $P$-orbit
projection  $f$, that is, $\kappa$ closes the following diagram up
to commutative: $\kappa\circ\varphi=\psi\circ f$.
  \end{lem}
 {\it Proof.} Let  $x=[g,s]_K$ and $x'=[g',s']_K
\in G\times_K \Bbb S$. If  $f(x)=f(x')$, then
$\pi(g')=\pi(g)\cdot l\ii$ and $f'(s')=l\cdot f'(s)$ for $l\in L$.

Since  $\pi'$ is the epimorphism, $l=\pi(k),k\in K$, and therefore
$f'(s')=\pi(k)\cdot f'(s)=f'(k\cdot s)$. Since $f'$ is the
$Q$-orbit projection, $s'$ and $k\cdot s$ lie  on the same
$K$-orbit. Hence it is easily  deduced that $x'$ and $x$ lie on
the same $G$-orbit.

All other characterization properties  for $P$-orbit projections
from Proposition  \ref{Proj1} is checked straightforwardly  and we
leave it to the reader. $\ \square$

\medskip  
The
following theorem  on {\it extension  of  tube  structure of
maps}
has a highly important role in inductive argument.

 \begin{thm}\label{thre+e+}
Let $\mathcal B$ be a commutative  diagram,
   $$
\begin{array}{ccc}
\Bbb X=\op{Cl} \Bbb X & \stackrel{j}{\hookrightarrow} & \Bbb W\\
f\downarrow&&\downarrow \hat f\\ \Bbb Y=\op{Cl} \Bbb
Y&\stackrel{i}{\hookrightarrow} &\Bbb Z
\end{array}
$$
 in which  $\mps{f}{\Bbb X}{\Bbb Y}$  and $\mps{\hat f}{\Bbb
W}{\Bbb Z}$ are  $P$-orbit  projections. If  $f$ has a
$\kappa$-tube structure generated  by slice maps
$\mps{\varphi}{\Bbb X}{G/K}$ and $\mps{\psi}{\Bbb Y}{H/L}$, then
there exist  an invariant neighborhoods  $\Bbb B,\Bbb Y\subset\Bbb
B\subset\Bbb Z$, and $\Bbb A\mean\hat f\ii(\Bbb B),\Bbb
X\subset\Bbb A\subset\Bbb W$, such that the $P$-orbit  projection
$\mps{f\restriction}{\Bbb A}{\Bbb B}$ has a $\kappa$-tube
structure generated  by slice maps $\mps{\hat\varphi}{\Bbb
A}{G/K}$ and $\mps{\hat\psi}{\Bbb B}{H/L}$, which  are  extensions
of $\varphi$ and $\psi$ respectively.
 \end{thm}
 {\it Proof  of Theorem \ref{thre+e+}.}
Let  $\mathcal A$ be a diagram which generates the $\kappa$-tube
structure on $f$. Since the subgroup  $H<G$ is extensor, it
follows by $(10)$ from Section 4 that $G$ admits an orthogonal
action on $\Bbb R^n$ such that $G/K$ is equimorphic to the orbit
of a point. Hence we can assume without loss of generality that
$G/K$ is  an orbit  in $H/L\times\Bbb R^N$, moreover, $\kappa$
coincides with the restriction of the projection
$\mps{\op{pr}_1}{H/L\times\Bbb R^N}{H/L}$ on $G/K$. Represent the
composition $\Bbb
X\matop{\rightarrow}{}{\varphi}G/K\hookrightarrow H/L\times\Bbb
R^N$ as $(\varphi_1,\varphi_2)$.

Since  by  $(2)$ the map  $\mps{\kappa}{G/K}{H/L}$ is
equivariantly locally soft, there exists  a fiberwise equivariant
retraction $\mps{r}{\Bbb V}{G/K}$ of an invariant neighborhood
$\Bbb V\subset H/L\times \Bbb R^N$ of $G/K$ such that
$\op{pr}_1\circ r=\op{pr}_1$. In view of $H/L\in H$-$\op{ANE}$ and
$\Bbb R^N\in G$-$\op{AE}$ there exist a local $H$-extension
$\mps{\tilde{\psi}}{\Bbb B'}{H/L,\Bbb Y\subset\Bbb B'\subset\Bbb
Z}$, of $\psi$ and a $G$-extension $\mps{\chi}{\Bbb W}{\Bbb R^N}$
of $\mps{\varphi_2}{\Bbb X}{\Bbb R^N}$.

Let   $\Bbb A\mean\hat f\ii(\Bbb B')$. Consider the $G$-map
$\mps{\sigma}{\Bbb A'}{H/L\times \Bbb R^N}$ given by the  formula
$\sigma=(\tilde{\psi}\circ\hat f )\times\chi$. It is clear that
$\Bbb A\mean\sigma\ii(\Bbb V)$ and $\Bbb B\mean \hat f(\Bbb A)$
are invariant neighborhoods of $\Bbb X$ and $\Bbb Y$ respectively.
We assert that $\mps{\hat\varphi\mean r\circ\sigma}{\Bbb A }{G/K}$
and $\hat\psi\mean\tilde{\psi}\restriction_{\Bbb B}$ generate a
$\kappa$-tube  structure on $\hat f\restriction_{\Bbb A }$:
$\kappa\circ \hat\varphi=\hat f\restriction\circ \hat \psi$. $ \ \
\square$

%%%%%%%%%%%%%%%%%%%%%%%%%%%%%%%%%%%%%%%%%%%%%%%%%%%%%%%%%%%%%
%%%%%%%%%%%%%%%%%%%%%%%%%%%%%%%%%%%%%%%%%%%\input{EqEmb4}

\section{Equivariant homotopy density}

 We postpone the proof  of Theorem \ref{v.03.04+3} to the end of the
section in view of the necessity of certain auxiliary facts. For
a 
compact group  $G$ we denote by  $\T$  the discrete union of all
transitive spaces $G/H\in G$-$\ANE$. It is clear that $\T$ is
metrizable, each point of the metric cone  $\op{Con} \T$ over $\T$
is
an
extensor and its orbit space  is the cone  over
a discrete space.
   \begin{thm} \label{v.03.04+}
 Let $\x$ be a metric $G$-subspace of the $G$-space $\Bbb Y$ whose orbit
space $Y$ is metrizable. Then for each nested family $\{\Bbb
V_n\subset\Bbb Y\}_{n=1}^\infty$ of invariant neighborhoods  of
$\x$ there exists a $G$-map $\mps{f}{\Bbb Y}{\T_G}\mean(\op{Con}
\T)^\omega$ such that
 \begin{enumerate}
 \item[$(a)$] the restriction of $f$ on $\Bbb X$ is
isovariant;
 \item[$(b)$] $f(\Bbb Y\setminus\cap\{\Bbb
V_n\mid n\ge 1\})\subset (\T_G)_{\mathcal E}$.
 \end{enumerate}
  \end{thm}
 {\it Proof  of Theorem \ref{v.03.04+}.}
We can assume without loss of generality that
  \begin{enumerate}
 \item[$(c)$] $X$ contains no isolated points.
 \end{enumerate}
 It otherwise should pass from $\Bbb X$ and $\Bbb Y$ to $\Bbb
X\times [0,1]^\omega$ and $\Bbb Y\times [0,1]^\omega$
respectively.

Since  $Y=\Bbb Y/G$ is metric, there exists  a family  ${\mathcal
B} =\{W_{\mu}\}_{\mu \in M}\ $ of open  subsets of $Y$
intersecting $X$ such that
 \begin{enumerate}
 \item ${\mathcal B}=\matop{\bigcup}{n=1}{\infty}{\mathcal B}_n$, where  ${\mathcal B}_n
=\{W_{\mu}\}_{\mu\in M_n\subset M}$ is a discrete family,
$\matop{\coprod}{n=1}{\infty}M_n=M$;
 \item the body of ${\mathcal B}_n$ is contained in $\Bbb
V_n$ for each  $n$;
 \item the restriction  ${\mathcal B}\restriction_{X}$ generates the basis of  $X$.
 \end{enumerate}

 By  $\Bbb W_\mu$ we denote  $p^{-1}W_{\mu}$ where
$\mps{p}{\Bbb Y}{Y}$ is the orbit  projection. Since $\Bbb W_\mu$
has a trivial slice, the number
 $$i(\mu)\mean\inf\{\diam
\big( \Bbb X\cap \varphi^{-1}(g\cdot[H])\big) \mid
\varphi\colon\Bbb W_\mu\to G/H\ \text{is a slice map  and }\ g\in
G\}\ge 0 $$
 is correctly  defined  (here we take the diameter  with respect to the compatible invariant
metric  $\varrho$ existing on $\Bbb X$ by \cite{Pal}). By $(c)$,
$i(\mu)>0$. It follows by invariance of $\sigma$ that the
diameters of $\Bbb X\cap \varphi^{-1}(g\cdot[H])=\Bbb X\cap
g\cdot\varphi^{-1}([H])$ and $\Bbb X\cap \varphi^{-1}([H])$ are
equal. Hence it is sufficient to take $g=e$ in the definition of
$i(\mu)$. Of particular interest is those slice map
$\varphi_\mu\colon\Bbb W_\mu\to G/H_\mu,\mu\in M_n$, for which
$\diam \big( \Bbb X\cap \varphi_\mu^{-1}([H_\mu])\big)<
j(\mu)\mean 2i(\mu)$. It is easy to see that $\varphi_\mu$
satisfies the following important property.
   \begin{lem}\label{eqemb4} If  $\mu\in M_n$, then
$\ \ \ \diam\big(\Bbb X\cap \varphi_\mu^{-1}([H_\mu])\big)<2\diam
\big(\Bbb X\cap \varphi^{-1}([H])\big)$ for each slice map
$\varphi\colon\Bbb W_\mu\to G/H$.
  \end{lem}
  We consider the $G$-map  $\mps{\psi_\mu\mean
\op{Con}\varphi_\mu} {\Bbb Y}{\op{Con}G/H_\mu}$ coinciding with
$(\varphi_\mu,\xi_\mu)$ on $\Bbb W_\mu$, and  with the vertex
$\{\ast\}$  on the complement to $\Bbb W_\mu$ (here $\xi_\mu:\Bbb
Y\to[0,1]$ is a function  constant on orbits and such that
$\xi_\mu^{-1}(0)=\Bbb Y\setminus \Bbb W_\mu$). It is clear that
$\psi_\mu\ii(G/H_\mu\times(0,1])=\Bbb W_\mu$ and
$G_{\varphi_{\mu}(y)}=G_{\psi_{\mu}(y)}$ for all  $y\in\Bbb
W_\mu$.

Since the family ${\mathcal B}_n$ is discrete, it is easy to see
that the formulae $\psi_{n}\restriction_{\Bbb W_{\mu}}= \psi_\mu$
for $\mu\in M_{n}$ and  $\psi_{n}\restriction_{\Bbb
Y\setminus\cup\{\Bbb W_{\mu}\mid \mu\in M_{n}\}}=\{\ast\}$
correctly  define the  continuous $G$-map  $\psi_{n}:\Bbb Y\to
\op{Con } \T$. It is clear  that $\psi_n\ii(*)=\Bbb
Y\setminus\cup\{\Bbb W_{\mu}\mid \mu\in M_{n}\}$. It turns out
that the desired $G$-map $f$ is the diagonal product
$\Delta\psi_{n}:\Bbb Y\to\big(\op{Con} \T\big)^\omega =\T_G$.

Fix arbitrary  $x\in\x$ and $\alpha>0$. By Theorem \ref{slice2-6}
on approximate slice there exists a $G$-map $r:\Bbb U(x)\to G/H\in
G$-$\ANE$ of some neighborhood  $\Bbb U(x)\subset\Bbb Y$ for which
$r^{-1}([H])$ has  diameter less then $\alpha/2$. Since  by $(3)$
the restriction ${\mathcal B}\restriction_{X}$ is the basis of
$X$, there exists an index  $\mu\in M_n$ such that $x\in \Bbb
W_{\mu}\subset \Bbb U(x)$. We note that $\diam
(r^{-1}([H]))<\alpha/2$ implies $i(\mu)<\alpha/2$, and hence it
follows by Lemma \ref{eqemb4} that
 \begin{enumerate}
 \item[$(d)$]
$\diam\big(\Bbb X\cap
\varphi_\mu^{-1}([H_\mu])\big)<i(\mu)<\alpha$.
 \end{enumerate}
 Therefore $\diam \big(\Bbb X\cap \varphi_\mu^{-1}(g\cdot[H_\mu])\big)<\alpha$
for all  $g\in G$. Hence it follows the existence of a sequence of
slice maps   $\{\varphi_{\mu_i}\colon\Bbb W_{\mu_i}\to
G/H_{\mu_i}\mid i\ge 1,\mu_i\in M_{n_i}\}$ such that for
$a_{i}\mean\varphi_{\mu_i}(x),i\ge 1$, we have
$\varphi_{\mu_i}^{-1}(a_{i})\subset\Bbb W_{\mu_i}$  and
 \begin{enumerate}
 \item[$(e)$] $\{\diam\big(\Bbb X\cap\varphi_{\mu_i}^{-1}(a_{i})\big)\}\rightarrow0$.
  \end{enumerate}
  Since $\varphi_{\mu_i}(x)=a_{i}$ and $G_{a_i}=G_{\varphi_{\mu_i}(x)}=G_{\psi_{\mu_i}(x)}\supset
G_{f(x)}\supset G_x$, $G_x\subset\cap\{G_{a_i}\mid i\ge 1\}$. To
prove that the stabilizers of $x$ and $f(x)$ are equal, it is
sufficient to establish the converse inclusion $\cap\{G_{a_i}\mid
i\ge 1\}\subset G_x$. If $g\in\cap\{G_{a_i}\mid i\ge 1\}$, then
$g\cdot x\in \varphi_{\mu_i}^{-1}(a_{i})$ for all $i\ge 1$, and by
virtue of $(e)$, we have $x=g\cdot x$. Hence $g\in G_x$.

 It now remains to check $(b)$. If  $y\not\in\Bbb
V_n$, then $y\not\in\Bbb V_{n+m}$ and hence it follows from $(2)$
that all coordinates of $f(y)$ excepting the first $n$ ones
coincide with the vertex $\{*\}$. Hence $G_{f(y)}$ is an extensor
subgroup as the intersection of finitely many of extensor
subgroups (Proposition \ref{aslice2-6}). $\ \ \square$

\bigskip Theorem  \ref{v.03.04+} implies an important  result
on the structure of solutions of the closed  PEA.

   \begin{prop} \label{v.03.04+5}
For each  solution  $s'\colon\Bbb X\hookrightarrow\Bbb
Y'\in{\mathcal M}$ of PEA  for closed  admissible  diagram
$\mathcal D$ there exists a solution  $s\colon\Bbb
X\hookrightarrow\Bbb Y$ of PEA for $\mathcal D$ such that $s\ge
s'$ and $\Bbb Y\setminus\Bbb X\subset{\Bbb Y}_{{\mathcal E}}$.
  \end{prop}
  
{\it Proof of Proposition  \ref{v.03.04+5}.} By  Theorem
\ref{v.03.04+} there exists  a $G$-map  $\mps{f}{\Bbb Y'}{\T_G}$
such that $f\restriction_{\Bbb X}$ is isovariant  and $f(\Bbb
Y'\setminus\x)\subset (\T_G)_{\mathcal E}$. We consider the
fiberwise  product  $\Bbb Y\mean Y_f\times_p \T_G\in{\mathcal M}$
where $\mps{p}{\T_G}{\T_G/G}$ is the orbit  projection.

Since  $f\restriction_{\Bbb X}$ is isovariant, $\Bbb X\subset\Bbb
Y$. It is easy to check that $s\colon\Bbb X\hookrightarrow\Bbb Y$
covers $X\hookrightarrow Y$ and the natural $G$-map $\mps{h}{\Bbb
Y'}{\Bbb Y},h(y')= (p(y'),f(y'))$, satisfies $h\circ s'=s$ and
$\tilde h=\op{Id}_{Y}$. It is clear  that $\Bbb Y\setminus\Bbb
X\subset{\Bbb Y}_{{\mathcal E}}$. $\ \square$

\medskip

We note that each closed $G$-embedding $s\colon\Bbb
A\hookrightarrow\Bbb Z$ is the solution of the closed PEA for
diagram ${\mathcal D}\mean\{\Bbb
A\matop{\rightarrow}{}{p}A\hookrightarrow Z\}$. As an easy
application of Proposition \ref{v.03.04+5} to ${\mathcal D}$, we
get
   \begin{prop} \label{v.03.04+7}
For each  partial   $G$-map $\pam{\mathcal M\ni\Bbb Z}{\Bbb
A}{f}{\Bbb X}\in G$-$\op{AE}(\mathcal M)$ there exists a
$G$-extension $\mps{\hat f}{\Bbb Z}{\Bbb X}$ such that $\hat
f(\Bbb Z\setminus\Bbb A)\subset{\Bbb X}_{{\mathcal P}}$. The
similar result  takes place for a $G$-$\op{ANE}$-space.
  \end{prop}
\bigskip To prove Theorem  \ref{v.03.04+3}, we
apply Proposition \ref{v.03.04+7} to the partial $G$-map
$\pam{\Bbb X\times[0,1]}{\Bbb X\times\{0\}}{\op{Id}}{\Bbb X}\in
G$-$\op{AE}$. The case $\Bbb X\in G$-$\op{ANE}$ is proved
analogically. $\ \ \square$

%%%%%%%%%%%%%%%%%%%%%%%%%%%%%%%%%%%%%%%%%%%%%%%%%%%%%%%%%%%%%%%%
%%%%%%%%%%%%%%%%%%%%%%%%%%%%%%%%%%%%%%%\input{ActExt22}

\section{Extensions of $P$-orbit projections}

  We treat the closed  problem of extending the
  action
in the general context. Let $\mps{f}{\x}{\Bbb Y}$ be a $P$-orbit
projection  for the kernel  $P$ of the epimorphism
$\mps{\pi}{G}{H}$ of compact groups, and $i\colon\Bbb
Y\hookrightarrow \Bbb Z$ an arbitrary $H$-embedding of $\Bbb Y$
into the equivariantly hereditarily paracompact, i.e.
each  open invariant  subset is paracompact (it is equivalent to
the hereditary paracompactness of the orbit space), $H$-space
$\Bbb Z$. The resulting diagram ${\mathcal D}\mean\{\Bbb
X\mathop{\rightarrow}\limits^{f}\Bbb
Y\mathop{\hookrightarrow}\limits^{i} \Bbb Z\}$ is called  {\it
$H$-admissible}. It is clear  that the induced diagram $\Bbb
X\mathop{\longrightarrow}\limits^{p_{\Bbb Y}\circ f}
Y\mathop{\hookrightarrow}\limits^{\tilde i}Z$ is $G$-admissible or
merely admissible. We say that
 \begin{enumerate}
 \item {\it The general problem on extending of action  (GPEA) is solvable
for $H$-admissible  diagram } ${\mathcal D}$ if  there exist a
$G$-embedding  $j:\Bbb X\hookrightarrow\Bbb W$ into a $G$-space
$\Bbb W$  and a $P$-orbit  projection  $\mps{\hat{f}}{\Bbb W}{\Bbb
Z}$ such that $\hat{f}\circ j=i\circ f$;
 \item The {\it GPEA is locally solvable
for $H$-admissible diagram} ${\mathcal D}$ if it is solvable for
some $H$-admissible diagram  $\Bbb
X\mathop{\rightarrow}\limits^{f}\Bbb Y
\mathop{\hookrightarrow}\limits^{i} \Bbb E$ where   $\Bbb Y\subset
\op{Int}\Bbb E\subset\Bbb Z$;
 \item The {\it general  problem
on extending  of action  is solvable  for $P$-orbit projections}
if GPEA is solvable  for each  $H$-admissible diagram.
 \end{enumerate}

 It will be shown below that the main result  --
Theorem  \ref{one}, is reduced to the following key theorem on
solvability of the closed GPEA.
 \begin{thm}\label{one++} The closed GPEA is solvable  for all $P$-orbit
projections, provided that the kernel  $P$ is a compact Lie group.
 \end{thm}
 We note that the content of Theorems  
 \ref{one} and \ref{one++} is identical in the case of $G=P$.
Fist we consider the simplest case of Theorem \ref{one++} and give
the complete its proof in the following section.
   \begin{lem}\label{three+} The GPEA for each $H$-admissible diagram  $\Bbb
X\mathop{\rightarrow}\limits^{f}\Bbb
Y\mathop{\hookrightarrow}\limits^{i} \Bbb Z$ is solvable, provided
$i$ is open.
 \end{lem}
  Hence it easily follows  that
  \begin{enumerate}\setcounter{enumi}{3}
 \itemsep=-2mm
 \item If  the GPEA is locally  solvable for $H$-admissible
diagram ${\mathcal D}$,  then the GPEA is solvable for this
diagram ${\mathcal D}$.
 \end{enumerate}
 {\it Proof of Lemma \ref{three+}.} We construct the $G$-space  
 $\Bbb W$ in a such
manner that  $\Bbb W\setminus\Bbb X=\Bbb Z\setminus \Bbb Y$. 
With this goal
in view, we set  $\Bbb W\mean\Bbb X\sqcup(\Bbb Z\setminus
\Bbb Y)$. 
The basis of $\Bbb W$ is generated by all open sets in
$\Bbb X$ and by sets  $\{\tilde O\mean f^{-1}(O\cap \Bbb
Y)\sqcup(O\setminus \Bbb Y)\mid O\subset \Bbb Z\text{is arbitrary
open set}\}$. The action of $G$ on $\Bbb W$ is defined as 
follows.
It coincides with the action of $G$ on $\Bbb X$, and it
is given by the formula  $g\cdot y=\pi(g)\cdot y,y\in \Bbb
Z\setminus \Bbb Y$, on $\Bbb Z\setminus \Bbb Y$. We can easily
verify the continuity of this action. It is clear that the $G$-map
$\mps{\hat{f}}{\Bbb W}{\Bbb Z}$ coincides  with $f$ on $\x$ and
with $\Id$ on $\Bbb Z\setminus \Bbb Y$ is a $P$-orbit projection
extending  $f$.  $\ \square$

 \bigskip
Let  $\Bbb X\mathop{\rightarrow}\limits^{f}\Bbb
Y\mathop{\hookrightarrow}\limits^{i} \Bbb Z$ be an $H$-admissible
diagram. We note that $\Bbb A=f\ii( f(\Bbb A))$ for all $\Bbb
A\subset\x$, and also  $G_x=P\cdot G_x=\pi\ii(\pi(G_x))$ for
$P<G_x$. Hence  it follows by perfectness of $f$ that
\begin{enumerate}\setcounter{enumi}{4}
 \item $\mps{f\restriction_{\x^P}}{\x^P}{f(\x^P)}$
is an  equimorphism of  closed subsets of  $\Bbb X$ and $\Bbb Y$.
\end{enumerate}

 Let $\Bbb X'\mean\Bbb X\setminus \Bbb X^P$, $\Bbb
Z'\mean\Bbb Z\setminus \op{Cl}_{\Bbb Z}(f(\Bbb X^P))$  and $\Bbb
Y'\mean \Bbb Y\cap\Bbb Z'=\Bbb Y\setminus f(\Bbb X^P)$. It is
clear  that
\begin{enumerate}\setcounter{enumi}{5}
 \item  $(\Bbb X')^P=\varnothing$ and the
map  $\mps{f\restriction}{\Bbb X'}{\Bbb Y'}$ is a $P$-orbit
projection.
\end{enumerate}

 The following  assertion reduces the investigation of the closed GPEA
to the case of absence of $P$-fixed points in $\x$.
 \begin{lem}\label{four+-} If the closed GPEA for all $H$-admissible
diagram $\Bbb X\mathop{\rightarrow}\limits^{f}\Bbb
Y\mathop{\hookrightarrow}\limits^{i} \Bbb Z$ with empty $P$-fix
point set $\Bbb X^P$ is solvable,  then it is solvable for all
$H$-admissible diagram.
 \end{lem}
 {\it Proof.}
Since  $Z'$ is hereditarily paracompact, the diagram $\Bbb
X'\mathop{\rightarrow}\limits^{f\restriction}\Bbb
Y'\mathop{\hookrightarrow}\limits^{i'} \Bbb Z'$ is $H$-admissible.
Since   $(\Bbb X')^P=\varnothing$, the lemma can be
applied and hence
the closed GPEA is solvable: there exist a $G$-embedding $j':\Bbb
X'\hookrightarrow\Bbb W'$ into $G$-space $\Bbb W'$ and a $P$-orbit
projection $\mps{\hat{f}}{\Bbb W'}{\Bbb Z'}$ such that
$\hat{f}\circ j'=i'\circ f\restriction$. We apply Lemma
\ref{three+} to the natural open $H$-embedding  $i'':\Bbb Z'
\hookrightarrow \Bbb Z$: the GPEA for diagram $\Bbb
W'\mathop{\rightarrow}\limits^{\hat f}\Bbb
Z'\mathop{\hookrightarrow}\limits^{i''} \Bbb Z$ is solvable.
Therefore the  GPEA is solvable for arbitrary $H$-admissible
diagram. $\ \ \square$

\bigskip Let  ${\mathcal
D}=\{\Bbb X\mathop{\rightarrow}\limits^{f}\Bbb
Y\mathop{\hookrightarrow}\limits^{i} \Bbb Z\}$ be a closed
$H$-admissible diagram. Each representation of $\Bbb Z$ as a union
of closed $G$-subspaces  $\Bbb Z_1$ and $\Bbb Z_2$  generates
three closed $H$-admissible diagrams $\Bbb
X_i\mathop{\rightarrow}\limits^{f\restriction}\Bbb
Y_i\mathop{\hookrightarrow}\limits^{i\restriction} \Bbb
Z_i,i=0,1,2$, where  $\Bbb Z_0\mean \Bbb Z_1\cap \Bbb Z_2$, $\Bbb
Y_i\mean\Bbb Z_i\cap\Bbb Y$ and $\Bbb X_i\mean f\ii(\Bbb Y_i)$.
Suppose that $s_i\colon \Bbb X_i\hookrightarrow \Bbb W_i,i=0,1,2$,
are solutions of the closed GPEA  for these $H$-admissible
diagrams, so that $s_0$ is the restriction of $s_j$ on $\Bbb X_0$
for each $j=1,2$ (assuming that $\Bbb W_0$ is naturally contained
in $\Bbb W_j$ ). The following fact is evident:
 \begin{lem}\label{four+-+} Let
$\Bbb W$ be a natural gluing of $G$-spaces $\Bbb W_1$ and $\Bbb
W_2$ along $\Bbb W_0$. Then $\Bbb W_i\subset\Bbb W,i=1,2$, are
closed invariant subsets of $G$-space $\Bbb W$, and there exists a
solution $s\colon\Bbb X\hookrightarrow\Bbb W$  of the closed GPEA
for $\mathcal D$ such that $s_i$ is the restriction of $s$ on
$\Bbb X_i$ for each $i=0,1,2$.
 \end{lem}
As an easy corollary  of Lemma \ref{four+-+} we get
 \begin{prop}\label{four+-+-} Let $\Bbb F$ be a closed $G$-subset of $\Bbb
Z$ such that $\Bbb Z=\Bbb Y\cup\Bbb F$. If the closed GPEA  is
solvable  for $H$-admissible diagram $f\ii(\Bbb X\cap \Bbb
F)\mathop{\rightarrow}\limits^{f\restriction}\Bbb Y\cap\Bbb
F\mathop{\hookrightarrow}\limits^{i} \Bbb F$, then the closed GPEA
is solvable  for ${\mathcal D}$.
 \end{prop}
 
 \medskip
 In conclusion, we explain the  
 {\bf reduction of Theorem \ref{one} to Theorem \ref{one++}.}
 
   \begin{prop}   \label{End1+*}
 The validity of Theorem  \ref{one++} implies the validity of Theorem \ref{one}.   \end{prop}
    {\it Proof.}
Let 
$\{ P_\alpha\triangleleft G\}$ 
be a Lie series of $G$ with
$P_1=G$. We represent $\x$ as $\matop{\lim}{\leftarrow}{} \{
\x_\alpha, f^{\beta}_{\alpha}\}$ 
(see the notation from Lemma \ref{End1}). Since $\x_1=X$ and $\Bbb Y_1=Y$, the embedding
$i\colon X\hookrightarrow Y$ can be identified with a
$G/P_1$-embedding $i_1\colon\Bbb X_1\hookrightarrow\Bbb Y_1$.
Since $P_\alpha/P_{\alpha+1}$ is a compact Lie group, Theorem
\ref{one++} can be applied many times. Hence it follows by
transfinite induction on $\alpha$ that there exist a
$G/P_\alpha$-embedding $i_\alpha:\Bbb X_\alpha\hookrightarrow\Bbb
Y_\alpha$ into $G/P_\alpha$-spaces $\Bbb Y_\alpha$  and
$P_{\alpha}^{\alpha+1}$-orbit  projection
$\mps{\hat{f}_{\alpha}^{\alpha+1}}{\Bbb Y_{\alpha+1}}{\Bbb
Y_{\alpha}}$ such that $\hat{f}_{\alpha}^{\alpha+1}\circ
i_{\alpha+1}=i_{\alpha}\circ f_{\alpha}^{\alpha+1}$ for each
$\alpha<\omega$. Then $\x=\matop{\lim}{\leftarrow}{}\{ \Bbb
X_\alpha, f_\alpha^\beta\}$ naturally lies in the $G$-space  $\Bbb
Y\mean\matop{\lim}{\leftarrow}{}\{ \Bbb Y_\alpha,
\hat{f}_\alpha^\beta\}$ and $\x\hookrightarrow\Bbb Y$ is the
desired $G$-embedding. $\ \ \square$

 \medskip
%%%%%%%%%%%%%%%%%%%%%%%%%%%%%%%%%%%%%%%%%%%%%%%%%%%%%%%%%%%%
%%%%%%%%%%%%%%%%%%%%%%%%%%%%%%%%%%%%%%%%\input{ActExt23}
\section{Proof  of Theorem  \ref{one++}}

   We use the argument based on
Palais metatheorem  \ref{t5.9++}. If  $|P|=1$ then $\pi$ is the
isomorphism that trivializes the situation under consideration. We
suppose now that for each  proper subgroup $Q<P$ the closed GPEA
for each $Q$-orbit  projection  is solvable and show that it is
solvable for each  $P$-orbit projection $\mps{f}{\x}{\Bbb Y}$. By
Lemma \ref{four+-} it is sufficient  to study a $G$-space $\Bbb X$
without $P$-fixed points, $\Bbb X^P=\varnothing$.

First we consider the case of the $P$-orbit projection $f$ having
a nontrivial tube structure.
  \begin{lem}\label{five-1+} If  $f$ has
a nontrivial tube structure,  then the closed GPEA is solvable for
each $H$-admissible  diagram  $\Bbb
X\mathop{\rightarrow}\limits^{f}\Bbb Y\mathop{\hookrightarrow}
\limits^{i}\Bbb Z$.
 \end{lem}
{\it Proof of Lemma \ref{five-1+}.} Consider the slice maps
$\mps{\varphi}{\Bbb X}{G/K}$ and $\mps{\psi}{\Bbb Y}{H/L}$ from
commutative diagram $\mathcal A$ which generate the nontrivial
tube structure on $f$. Since $H/L\in G$-$\op{ANE}$, $\psi$ has an
$H$-extension $\mps{\hat\psi}{\Bbb U}{H/L}$ in some invariant
neighborhood $\Bbb U,\Bbb Y \subset\Bbb U\subset\Bbb Z$.

We note that by  Lemma \ref{Proj2} the map
$\mps{f\restriction}{\Bbb S=\varphi\ii[K]}{\Bbb T=\psi\ii[L]}$ is
a $Q$-orbit  projection for proper compact group  $Q<P$. Since
$\Bbb T'\mean\hat\psi\ii[L]$ is equivariantly hereditarily
 paracompact, the diagram $\Bbb
S\matop{\rightarrow}{}{f}\Bbb T\hookrightarrow\Bbb T'$ is
$L$-admissible. In view of the inductive hypothesis, the closed
GPEA is solvable for all $Q$-orbit projections, and hence it is
solvable for $\Bbb S\matop{\rightarrow}{}{f}\Bbb
T\hookrightarrow\Bbb T'$. Hence there exist a $K$-embedding
$j':\Bbb S\hookrightarrow\Bbb S'$ into $K$-space $\Bbb S'$ and a
$Q$-orbit  projection $\mps{f'}{\Bbb S'}{\Bbb T'}$ such that
$f'\circ j'=i\circ f\restriction$.

Next, by  Lemma \ref{Proj3} the map $\mps{\hat{f}=\pi\times
f'}{G\times_K \Bbb S'}{H\times_L\Bbb T'=\Bbb U}$ is a $P$-orbit
projection, and, as is easy to see, $\hat{f}$ solves the closed
GPEA for $H$-admissible  diagram  $\Bbb
X\mathop{\rightarrow}\limits^{f}\Bbb Y
\mathop{\hookrightarrow}\limits^{i} \Bbb U$. We apply  Lemma
\ref{three+} for diagram $G\times_K \Bbb
S'\mathop{\rightarrow}\limits^{\hat f}H\times_L\Bbb T'=\Bbb
U\mathop{\hookrightarrow} \Bbb Z$ with open  embedding $\Bbb
U\mathop{\hookrightarrow} \Bbb Z$ and complete the proof. $\ \
\square$

 \bigskip The last case of the proof of
Theorem  \ref{one++} consists in consideration of a $P$-orbit
projection $\mps{f}{\x}{\Bbb Y}$ with $\Bbb X^P=\varnothing$ and a
closed $H$-embedding  $\Bbb Y \mathop{\hookrightarrow}\limits^{i}
\Bbb Z$.
   \begin{lem}\label{ActExt3++}
There exist a closed neighborhood $\Bbb E,\Bbb Y\subset \Bbb
E\subset \Bbb Z$ and  its locally finite closed $H$-cover
$\sigma=\{\Bbb F_\gamma\subset \Bbb
E\}_{\gamma\in\Gamma}\in\cov\Bbb E$ such that for each
$\gamma\in\Gamma$
 \begin{enumerate}
 \itemsep=-2mm
\item the map  $g_\gamma \mean f\restriction\colon\Bbb
V_\gamma\matop{\longrightarrow}{}{}\Bbb F_\gamma\cap \Bbb Y$ has a
nontrivial  tube  structure where  $\Bbb V_\gamma\mean f\ii(\Bbb
Y\cap \Bbb F_\gamma)$.
\end{enumerate}
 \end{lem}
 {\it Proof.} Since  $\Bbb X^P=\varnothing$,
Proposition  \ref{Proj2-+} implies that for each $x\in\x$,
$f\restriction_{G(x)}$  has a nontrivial  tube structure. Hence it
follows by Theorem \ref{thre+e+} on extension  of  tube  structure
of maps that there exists a locally finite open $H$-cover
$\omega=\{\Bbb U_\alpha \}\in\cov\Bbb Y$ such that
\begin{enumerate}\setcounter{enumi}{1}
  \item  Each $P$-orbit  projection
$\mps{f\restriction}{f\ii(\Bbb U_\alpha)}{\Bbb U_\alpha}$ has a
nontrivial  tube  structure. \end{enumerate}

Let  $\nu$ be a family of open  $H$-sets of $\Bbb Z$ the
restriction of which on $\Bbb Y$ coincides  with $\omega$. Since
$Z$ is hereditarily paracompact, the body $\cup\nu$ is
paracompact. Hence  there exists a closed locally finite $H$-cover
$\sigma'$ of $\cup\nu$ refining $\nu$. It is clear that
$\sigma\mean\{\Bbb F\in\sigma'\mid \Bbb F\cap \Bbb
Y\not=\emptyset\}$ and $\Bbb E\mean\cup\sigma$ are desired. $\
\square$

\medskip To argue by the new transfinite induction, we well--order
the set $\Gamma$ indexing the elements of $\{\Bbb F_\gamma\}$.
Without loss of generality we can assume that $\Gamma$ has the
maximal element  $\omega$. If $\Bbb Q_\gamma\mean \cup\{ \Bbb
F_{\gamma'}|\gamma'\le\gamma\}$, then it is obvious that $\Bbb
Q_\omega=\Bbb E$ is the body of locally finite increasing system
of closed  subsets $\{\Bbb Q_\gamma\}$, moreover  $\Bbb
Q_{\gamma'}\cup \Bbb F_{\gamma}=\Bbb Q_\gamma$ for each
$\gamma=\gamma'+1$.

By  transfinite induction on $\gamma$ we construct closed
neighborhoods  $\Bbb R_{\gamma},\Bbb Y\subset \Bbb R_\gamma\subset
\Bbb Q_\gamma\cup \Bbb Y$, such that
\begin{enumerate}\setcounter{enumi}{1}
  \item If  $\gamma_0$
is the  minimal element of $\Gamma$, then  $\Bbb R_{\gamma_0}=\Bbb
Y\cup \Bbb F_{\gamma_0}$;
 \item   For all  $\gamma\in\Gamma$, the closed GPEA for $H$-admissible  diagram  $\Bbb
X\mathop{\rightarrow}\limits^{f}\Bbb Y\mathop{\hookrightarrow}
\limits^{i_{\gamma}}\Bbb R_{\gamma}$ is solvable, that is, there
exist a closed $G$-embedding
$j_{\gamma}\colon\x\hookrightarrow\Bbb W_{\gamma}$ into a
$G$-space $\Bbb W_{\gamma}$ and a $P$-orbit projection  $\mps{\hat
f_{\gamma}}{\Bbb W_{\gamma}}{\Bbb R_{\gamma}}$ such that $\hat
f_{\gamma}\circ j_{\gamma}=i\restriction\circ f$;
    \item $\Bbb R_{\gamma_1}\subset \Bbb R_{\gamma_2}$,
$\Bbb W_{\gamma_1}\subset \Bbb W_{\gamma_2}$ and $\hat
f_{\gamma_2}\restriction_{\Bbb W_{\gamma_1}}=\hat f_{\gamma_1}$
for all  $\gamma_1<\gamma_2$;
  \item $\Bbb R_{\gamma_2}\setminus\Bbb R_{\gamma_1}\subset
  \cup\{\Bbb F_{\gamma}\mid \gamma_1<\gamma\le\gamma_2\}$ for all  $\gamma_1<\gamma_2$.
  \end{enumerate}

 By $(5)$ the  {\sl stabilization condition} of
 constructed neighborhoods $\Bbb R_{\gamma}$ follows: if
$\Gamma_{\Bbb A}\mean\{\gamma\in\Gamma\mid \Bbb A\cap \Bbb
F_{\gamma}\not=\emptyset\},\Bbb A\subset \Bbb Z$, is finite, then
for each  $\gamma'\in\Gamma$ the intersection $\Bbb A\cap \Bbb
R_{\gamma'}$ coincides  with $\Bbb A\cap\Bbb R_{\gamma_m}$ where
$\gamma_m\mean\max\{\gamma\in\Gamma_\Bbb A\mid
\gamma\le\gamma'\}$.

\medskip
It is clear  that $\Bbb R_{\omega}\subset \Bbb Q_\omega=\Bbb E$ is
a closed neighborhood of $\Bbb Y$ in $\Bbb Z$. Hence it follows by
$(3)$ that the closed GPEA for diagram ${\mathcal D}=\{\Bbb
X\mathop{\rightarrow}\limits^{f}\Bbb Y\mathop{\hookrightarrow}
\limits^{i}\Bbb Z\}$ is locally solvable. We take into account the
property $(4)$ from the previous section and conclude that the
closed GPEA for ${\mathcal D}$ is solvable, which completes the {\sl
proof of Theorem \ref{one++}}.

 \medskip The {\sl
base} of the inductive argument is easily established with help of
Lemma \ref{five-1+}.
 \begin{lem}\label{five+} The closed GPEA for diagram  $\Bbb
X\mathop{\rightarrow}\limits^{f}\Bbb Y\mathop{\hookrightarrow}
\limits^{i_{\gamma_0}}\Bbb R_{\gamma_0}=\Bbb Y\cup \Bbb
F_{\gamma_0}$ is solvable.
 \end{lem}
 {\it Proof.} By
Lemmas  \ref{five-1+} and \ref{ActExt3++} the closed GPEA for
$H$-admissible diagram  $\Bbb V_{\gamma_0}
\mathop{\longrightarrow}\limits^{g_{\gamma_0}}\Bbb
F_{\gamma_0}\cap\Bbb Y\hookrightarrow \Bbb F_{\gamma_0}$ is
solvable. To complete the proof we apply Proposition
\ref{four+-+-} to $H$-admissible diagram $\Bbb
X\mathop{\rightarrow}\limits^{f}\Bbb
Y\mathop{\hookrightarrow}\limits^{i} \Bbb Y\cup\Bbb F_{\gamma_0}$.
$\ \square$

 \medskip The {\bf inductive step} consists of
 the following proposition.
\begin{lem}\label{five-+1}
Assume that for $\gamma'\in\Gamma,\gamma'<\omega$, we have defined
a closed neighborhood $\Bbb R_{\gamma'},\Bbb Y\subset \Bbb
R_{\gamma'}\subset \Bbb Q_{\gamma'}\cup \Bbb Y$, a $G$-embedding
$j_{\gamma'}\colon\Bbb X\hookrightarrow\Bbb W_{\gamma'}$ into a
$G$-space  $\Bbb W_{\gamma'}$  and a $P$-orbit  projection
$\mps{\hat{f}_{\gamma'}}{\Bbb W_{\gamma'}}{\Bbb R_{\gamma'}}$
satisfying  $(3)$. Then  for $\gamma=\gamma'+1\in\Gamma$
   \begin{enumerate}\setcounter{enumi}{5}
   \item There exist a  closed
neighborhood $\Bbb R_{\gamma},\Bbb Y\subset \Bbb R_{\gamma}\subset
\Bbb Q_{\gamma}\cup \Bbb Y$, a $G$-embedding $j_{\gamma}:\Bbb
X\hookrightarrow\Bbb W_{\gamma}$ into a $G$-space $\Bbb
W_{\gamma}$ and a $P$-orbit  projection
$\mps{\hat{f}_{\gamma}}{\Bbb W_{\gamma}}{\Bbb R_{\gamma}}$
satisfying  $(3)$ such that $\Bbb R_{\gamma'}\subset\Bbb
R_{\gamma}$, $\Bbb W_{\gamma'}$ naturally lies in $\Bbb
W_{\gamma}$, $\hat{f}_{\gamma}\restriction_{\Bbb
W_{\gamma'}}=\hat{f}_{\gamma'}$ and  $\Bbb R_{\gamma}\setminus\Bbb
R_{\gamma'}\subset\Bbb F_{\gamma}$.
\end{enumerate}
   \end{lem}
  {\it Proof of Lemma \ref{five-+1}.}
Since  $\mps{g_{\gamma}}{\Bbb V_{\gamma}}{\Bbb F_{\gamma}\cap\Bbb
Y}$ has a nontrivial  tube  structure, it follows by Theorem
\ref{thre+e+} that
\begin{enumerate}\setcounter{enumi}{6}
 \item There exists a closed invariant neighborhood
$\Bbb S\subset\Bbb R_{\gamma'}$ of $\Bbb F_{\gamma}\cap\Bbb
R_{\gamma'}$ such that the $P$-orbit projection $\mps{h\mean\hat
f_{\gamma'}\restriction_{\Bbb T }}{\Bbb T}{\Bbb S}$ where  $\Bbb
T\mean(\hat f_{\gamma'}) \ii(\Bbb S)\subset\Bbb W_{\gamma'}$ has a
nontrivial  tube  structure.
\end{enumerate}
Let  $\hat {\Bbb S}\subset\Bbb R_{\gamma'}\cup\Bbb F_{\gamma}$ be
a closed  neighborhood of $\Bbb F_{\gamma}\cap\Bbb R_{\gamma'}$
such that $\hat {\Bbb S}\cap\Bbb R_{\gamma'}=\Bbb S$. Hence it
follows by Lemma \ref{five-1+} that the closed GPEA is solvable
for $H$-admissible  diagram $\Bbb
T\mathop{\rightarrow}\limits^{h}\Bbb S\mathop{\hookrightarrow}
\limits^{}\hat {\Bbb S}$. To complete the proof we apply
Proposition \ref{four+-+-} to $H$-admissible diagram $\Bbb
W_{\gamma'}\mathop{\rightarrow}\limits^{\hat f_{\gamma'}}\Bbb
R_{\gamma'}\mathop{\hookrightarrow}\limits^{} \Bbb
R_{\gamma}\mean\Bbb R_{\gamma'}\cup\hat {\Bbb S}$. The
verification of all details is evident and we leave it to the
reader. $\ \square$

   \bigskip
  Let  $\gamma\in\Gamma$ be a limit ordinal.
Since  $\{\Bbb F_{\gamma'}\}_{\gamma'<\gamma}$ is a locally finite
$H$-cover of $\Bbb Y\cup\Bbb Q_\gamma$, it follows by $(5)$ that
the increasing system of $G$-subspaces $\{\Bbb
W_{\gamma'}\subset\Bbb W_{\gamma''}\}_{\gamma'\le\gamma''<\gamma}$
and  the system of $P$-orbit projections $\{\mps{\hat
f_{\gamma'}}{\Bbb W_{\gamma'}}{\Bbb
R_{\gamma'}}\}_{\gamma'<\gamma} $ are locally stabilized. In view
of this remark, we set $\Bbb R_\gamma\mean \cup\{\Bbb
R_{\gamma'}\mid\gamma'<\gamma\}$, $\Bbb W_\gamma\mean \cup\{\Bbb
W_{\gamma'}\mid\gamma'<\gamma\}$, and $\mps{\hat f_\gamma }{\Bbb
W_\gamma}{\Bbb R_{\gamma}}\mean\hat f_{\gamma'} $ on $\Bbb
W_{\gamma'}$ for all $\gamma'<\gamma$. It is easily checked that
$\Bbb X\hookrightarrow\Bbb W_\gamma$ and $\hat f_\gamma$ solve the
closed GPEA for $H$-admissible  diagram $\Bbb
X\mathop{\rightarrow}\limits^{f}\Bbb Y\hookrightarrow \Bbb
R_{\gamma}$ in such a manner that $\Bbb W_{\gamma'}\subset\Bbb
W_\gamma$ and $\hat f_{\gamma'}=\hat f_{\gamma}\restriction_{\Bbb
W_{\gamma'}}$ for all $\gamma'<\gamma$.

\medskip
%%%%%%%%%%%%%%%%%%%%%%%%%%%%%%%%%%%%%%%%%%%%%%%%%%%%%%%%%%%
%%%%%%%%%%%%%%%%%%%%%%%%%%%%%%%%%%%%%%%%%%\input{ActExt24}
\section{  Proof  of Theorem  \ref{one*}. }
  By  Theorem  \ref{one} there exists a solution  $s'\colon\Bbb
X\hookrightarrow\Bbb Y'$ covering  $X\hookrightarrow Y$. Since $Y$
is paracompact and $\mps{p_{\Bbb Y'}}{\Bbb Y'}{Y}$ is perfect,
$\Bbb Y'$ is also paracompact.

Since  $\Bbb X$ is  stratifiable, there exists a continuous
bijection of $\Bbb X$ onto a metric space  \cite{Bor}. Hence there
exists a continuous map $\mps{\varphi}{\Bbb X}{B}$ to the Banach
space $B$ such that
\begin{enumerate}\setcounter{enumi}{0}
 \itemsep=-2mm
 \item its  restriction  on each orbit is an  embedding.
 \end{enumerate}
Consider the map  $\mps{\psi}{\Bbb X}{\op C(G,B)}$ given by the
formula $\psi(x)(g)=\varphi(g\ii x), g\in G,x\in\Bbb X$, and which
is by $(1)$ isovariant. Since  by  Theorem \ref{v.1.4} $\op
C(G,B)$ is  $G$-$\ANE$ for the class of paracompact spaces and
$\Bbb Y'$ is paracompact, there exists a $G$-map
$\mps{\hat\psi}{\Bbb Y'}{\op C(G,B)}$ extending  $\psi$.

 We consider the fiberwise product  $\Bbb Y\mean Y_{\theta}\times_{p}
\op C(G,B)$ where  $\mps{p}{\op C(G,B)}{\op C(G,B)/G}$ is the
orbit projection and $\mps{\theta}{Y}{\op C(G,B)/G}$ is the map
induced by $\theta$. The $G$-space  $\Bbb Y$ is stratifiable as a
subset  of the product of two stratifiable spaces. Since $\psi$ is
isovariant, the $G$-space  $\Bbb X$ is naturally contained in the
$G$-space $\Bbb Y$, $s\colon\Bbb X\hookrightarrow\Bbb Y$, moreover
$s$ covers $X\hookrightarrow Y$. Finally, the natural $G$-map
$\mps{h}{\Bbb Y'}{\Bbb Y},h(y')=\big(p_{\Bbb
Y'}(y'),\hat\psi(y')\big)$, satisfies the properties $h\circ s'=s$
and $\tilde h=\op{Id}_{Y}$. $\ \square$

%%%%%%%%%%%%%%%%%%%%%%%%%%%%%%%%%%%%%%%%%%%%%%%%%%%%%%%%%
%%%%%%%%%%%%%%%%%%%%%%%%%%%%%%%%%%%%%%%\input{Liter}


\begin{thebibliography}{99}

\bibitem{Ab}
  H. Abels,
{\it Universal proper G-spaces,} 
Math. Z. {\bf 159} (1978), 143--158.

\bibitem{AgU}
 S. M. Ageev,
 {\it  An equivariant Dugundji theorem,} 
 Russian Math. Surveys {\bf 45}:5 (1990), 219--220.


\bibitem{Ag4}
S. M. Ageev,
 {\it   On the extension of action,} 
 Moscow Univ. Math. Bull.\ {\bf 47}:5 (1992),  17--19.
 
\bibitem{AgIzv}
 S. M. Ageev,
 {\it  Classification of $G$-spaces,} 
 Russian Acad. Sci. Izv. Math. {\bf 41}:3 (1993), 581--591.
 
 \bibitem{Ag11}
  S. M. Ageev,
  {\it  The problem of  extending of action and extensor properties of orbit spaces,}
  Talk at the Reseach Seminar on General Topology, Moscow State University, 18 March 1993.
  
 \bibitem{Ag1}
  S. M. Ageev,
  {\it  Extensor properties of orbit spaces and the problem of  continuation of an action,} 
  Moscow Univ. Math. Bull. {\bf 49}:1 (1994), 11--16.

 \bibitem{Ag2}
  S. M. Ageev,
  {\it  On a problem of Zambakhidze and Smirnov,}
  Math. Notes {\bf 58}:1--2 (1995),  679--684.

 \bibitem{AgRep}
  S. M. Ageev and D. Repov\v{s},
  {\it  The Jaworowski method in the problem of the preservation of extensor properties by an orbit functor,}
 Math. Notes\ {\bf 71}:3 (2002), 470--473.
 
\bibitem{AP}
 P. S. Alexandrov and B. A. Pasynkov,
{\it An Introduction to Dimension Theory} [in Russian], 
Nauka, Moscow, 1975. 

\bibitem{Bor}
C. J. R. Borges,  {\it  On stratifiable spaces,} Pacif. J. Math.\
{\bf 17}:1 (1966), 1--16.

\bibitem{Br}
 G. Bredon,
 {\it Introduction to Compact Transformation Groups},
 Pure and
Appl. Math. {\bf 46}, 
Academic Press, New York, London, 1972.

\bibitem{Cau}
R. Cauty,
 {\it  Une generalization du theoreme de Borsuk-Whitehead-Hanner aux espaces stratifiables,} 
 C. R. Acad.Sc.Paris {\bf 275}:4 (1972), 271--274.

\bibitem{Dick}
 T. tom Dieck, 
 {\it Transformation Groups and Representation Theory},
 Springer-Verlag, Berlin, 1979.

\bibitem{Dk}
C. H. Dowker,
 {\it  On a theorem of Hanner,}
 Ark. Mat. {\bf 2} (1952),  307--313.

\bibitem{Matu}
T. Matumoto,
 {\it Equivariant $K$-theory and Fredholm operator,} J.
Fac. Sci. Univ. Tokyo Sect. IA Math.  {\bf 18} (1971), 109--125.

\bibitem{Ma}
 C. Maxwell,
 {\it  Homeomorphisms of topological trabsformation groups  into function spaces,} 
 Duke Math. J. {\bf 33} (1966), 567--574.
 
\bibitem{Mu}
M. Murayama,
 {\it On  G-ANRs and their G-homotopy types,}
 Osaka J. Math.  {\bf 20}:3 (1983), 479--512.
 
 \bibitem{Pal}
R. Palais,
 {\it   The classification of G-spaces,}
Memoirs of Amer. Math. Soc. {\bf  36} (1960).

\bibitem{Ptr}
L. S. Pontryagin, {\it Topological Groups},
Gordon and Breach, New York, 1966.

\bibitem{Sz}
 J. Szenthe,
 {\it On the topological characterization of transitive Lie group actions,}
 Acta Sci. Math. Szeged {\bf 36}:3-4 (1974), 323--344.

\bibitem{Va}
 R.  Vaught,
 {\it  On extending actions,}
Proc. Amer. Math. Soc. {\bf 107}:4 (1989), 1087--1090.

\bibitem{Vil}
  G. Villalobos,
 {\it  Inversion of the theorem on the existence of slices,}
Moscow Univ. Math. Bull. {\bf 55}:1 (2000), 26--29.


\end{thebibliography}
 \end{document}